\documentclass[10pt]{amsart}
\usepackage{amssymb}
\usepackage{latexsym}
\usepackage{euscript}

   \def\sE{{\mathfrak E}}   
   \def\sH{{\mathfrak H}}

\def\dA{{\mathbb A}}      \def\dC{{\mathbb C}}
\def\dD{{\mathbb D}}

   \def\dN{{\mathbb N}}   
      \def\dR{{\mathbb R}}

\DeclareMathOperator{\ran}{ran}          % range
\DeclareMathOperator{\dom}{dom}          % domain

\def\wt#1{{{\widetilde #1} }}

\def\wh#1{{{\widehat #1} }}

\def\bm\chi{\mbox{\boldmath$\chi$}}
\def\half{{\frac{1}{2}}}

\def\IM{{\rm Im\,}}
\def\ran{{\rm ran\,}}
\def\cran{{\rm \overline{ran}\,}}

\def\cdom{{\rm \overline{dom}\,}}

\def\dim{{\rm dim\,}}

\def\e{{\rm e}}

\def\cmr{{\dC \setminus \dR}}

\def\uphar{{\upharpoonright}}

\newtheorem{theorem}{Theorem}[section]
\newtheorem{proposition}[theorem]{Proposition}

\newtheorem{lemma}[theorem]{Lemma}

\theoremstyle{remark}
\newtheorem{remark}[theorem]{Remark}

\theoremstyle{definition}

\newtheorem{definition}[theorem]{Definition}

\numberwithin{equation}{section}

\begin{document}

\title[A general realization theorem]
{A general realization theorem for matrix-valued
Herglotz-Nevanlinna functions}

\author[S.V.~Belyi]{Sergey Belyi}
\address{Department of Mathematics\\
Troy State University\\
Troy, AL 36082, USA} \email{sbelyi@trojan.troyst.edu}
%\indent{\it URL:} http://spectrum.troyst.edu/$\sim$belyi}

\author[S.~Hassi]{Seppo Hassi}
\address{Department of Mathematics and Statistics \\
University of Vaasa \\
P.O. Box 700, 65101 Vaasa \\
Finland} \email{sha@uwasa.fi}

\author[H.S.V. de Snoo]{Henk de Snoo}
\address{Department of Mathematics and Computing Science\\
University of Groningen \\
P.O. Box 800, 9700 AV Groningen \\
Nederland} \email{desnoo@math.rug.nl}

\author[E.R.~Tsekanovski\u{\i}]{Eduard Tsekanovskii}
\address{Department of Mathematics\\
Niagara University, NY 14109 \\
USA} \email{tsekanov@niagara.edu}
%\curaddr{}
\dedicatory{Dedicated to Damir Arov \ on the occasion of his 70th
birthday \\ and to Yury Berezanski\u i on the occasion of his 80th
birthday}
\date{\today}%{July 23, 2005}
%\thanks{}
%\translator{}



\begin{abstract}
New special types of stationary conservative impedance and
scattering systems, the so-called non-canonical systems, involving
triplets of Hilbert spaces and projection operators, are
considered. It is
established that every matrix-valued Herglotz-Nevanlinna function
of the form
\[
 V(z)=Q+Lz+\int_{\dR} \left(
\frac{1}{t-z}-\frac{t}{1+t^2}\right)\, d\Sigma(t)
\]
can be realized as a transfer function of such a new type of
conservative impedance system. In this case it is shown that the
realization can be chosen such that the main and the projection
operators of the realizing system satisfy a certain commutativity
condition if and only if $L=0$. It is also shown that $V(z)$ with
an additional condition (namely, $L$ is invertible or $L=0$), can
be realized as a linear fractional transformation of the transfer
function of a non-canonical scattering $F_+$-system. In
particular, this means that every scalar Herglotz-Nevanlinna
function can be realized in the above sense.

Moreover, the classical Liv\v{s}ic systems
(Brodski\u{\i}-Liv\v{s}ic operator colligations) can be derived
from $F_+$-systems as a special case when $F_+=I$ and the spectral
measure $d\Sigma(t)$ is compactly supported. The realization
theorems proved in this paper are strongly connected with, and
complement the recent results by Ball and Staffans.
\end{abstract}

\maketitle

\section{Introduction}

An operator-valued function $V(z)$ acting on a finite-dimensional
Hilbert space $\sE$ belongs to the class of matrix-valued
Herglotz-Nevanlinna functions if it is holomorphic on $\cmr$,
if it is symmetric with respect to the real axis, i.e.,
$V(z)^*=V(\bar{z})$, $z\in \cmr$, and if it satisfies the
positivity condition
\[
 \IM V(z)\geq 0, \quad z\in \dC_+.
\]
It is well known (see e.g. \cite{GT}) that matrix-valued
Herglotz-Nevanlinna functions admit the following integral
representation:
\begin{equation}
\label{nev0}
 V(z)=Q+Lz+\int_{\dR}
      \left( \frac{1}{t-z}-\frac{t}{1+t^2}\right)\, d\Sigma(t),
\quad
    z \in \cmr,
\end{equation}
where $Q=Q^*$, $L\geq 0$, and $\Sigma(t)$ is a nondecreasing
matrix-valued function on $\dR$ with values in the class of
nonnegative matrices in $\sE$ such that
\begin{equation}
\label{int0}
 \int_{\dR} \frac{\left( d\Sigma(t)x,x \right)}{1+t^2} <\infty,
 \quad x \in \sE.
\end{equation}
The problem considered in this paper is the general operator
representation of these functions with an interpretation in system
theory, i.e., in terms of linear stationary conservative dynamical
systems. This involves new types of stationary conservative
impedance and scattering systems (non-canonical systems) involving
triplets of Hilbert spaces and projection operators. The exact
definition of both types of non-canonical systems is given below. It
turns out that every matrix-valued Herglotz-Nevanlinna function can
be realized as a matrix-valued transfer function of this new type of
conservative impedance system. Moreover, assuming an additional
condition on the matrix $L$ in \eqref{nev0} ($L$ is invertible or
$L=0$), it is shown that such a function is realizable as a linear
fractional transformation of the transfer matrix-valued function of
a conservative stationary scattering $F_+$-system. In this case the
main operator of the impedance system is the ``real part'' of the
main operator of the scattering $F_+$-system. In particular, it
follows that every scalar Herglotz-Nevanlinna function can be
realized in the above mentioned sense. This gives a complete
solution of the realization problems announced in \textit{``Unsolved
problems in mathematical systems and control theory''} \cite{HST4}
in the framework of modified Brodski\u{\i}-Liv\v{s}ic operator
colligations (in the scalar case via impedance and scattering
systems, in the matrix-valued case via impedance systems).
Furthermore, the classical canonical systems of the M.S.~Liv\v{s}ic type
(Brodski\u{\i}-Liv\v{s}ic operator colligations) can be derived from
$F_+$-systems as a special case when $F_+=I$ and the spectral
measure $d\Sigma(t)$ is compactly supported.

Realizations of different classes of holomorphic matrix-valued
functions in the open right half-plane, unit circle, and upper
half-plane play an important role in the spectral analysis of
non-self-adjoint operators, interpolation problems, and system
theory; see \cite{AlTs1}--\cite{BHST}, \cite{Br}--\cite{LiYa},
\cite{PoW}--\cite{WE}. For special classes of Herglotz-Nevanlinna
functions such operator realizations are known.

Consider, for instance, a matrix-valued Herglotz-Nevanlinna
function of the form
\begin{equation}
\label{nev1}
 V(z)=\int_{a}^{b}\frac{d\Sigma(t)}{t-z},
\quad z\in \cmr,
\end{equation}
with $\Sigma(t)$ a nondecreasing matrix-valued function on the
finite interval $(a,b)\subset \dR$. Then $V(z)$ has an operator
realization of the form
\begin{equation}
\label{real1}
 V(z)=K^*(A-zI)^{-1}K,
\quad z\in \cmr,
\end{equation}
where $A$ is a bounded self-adjoint operator acting on a Hilbert
space $\sH$ and $K$ is a bounded invertible operator from the
Hilbert space $\sE$ into $\sH$. Such realizations are due to
M.S.~Brodski\u{\i} and M.S.~Liv\v{s}ic; they have been used in the
theory of characteristic operator-valued functions as well as in
system theory in the following sense (cf. \cite{Lv1}-\cite{LvPo},
\cite{Br}, \cite{BrLv}, \cite{LiYa}). Let $J$ be a bounded,
self-adjoint, and unitary operator in $\sE$ which satisfies $\IM
A=KJK^*$. Then the aggregate
\begin{equation}
\label{BL}
 \Theta =\begin{pmatrix}
          A    &K  &J  \\
          \sH  &   &\sE
         \end{pmatrix}
\end{equation}
or
\begin{equation}\label{LivSys}
    \left\{%
\begin{array}{ll}
    (A-zI)x=KJ\varphi_-, & \hbox{} \\
    \varphi_+=\varphi_- -2iK^* x, & \hbox{} \\
\end{array}%
\right.
\end{equation}
is the corresponding, so-called, {\em canonical system} or {\em
Brodski\u{\i}-Liv\v{s}ic  operator colligation}, where $\varphi_-
\in \sE$ is an input vector, $\varphi_+ \in \sE$ is an output
vector, and $x$ is a state space vector in $\sH$. The function
$W_\Theta(z)$, defined by
\begin{equation}
\label{een+}
 W_\Theta(z)=I-2iK^*(A-zI)^{-1}KJ,
\end{equation}
such that $\varphi_+=W_\Theta(z) \varphi_-$, is the {\em transfer
function} of the system $\Theta$ or the {\em characteristic
function} of  operator colligation.  Such type of systems appear
in the theory of electrical circuits and have been introduced by
M.S.~Liv\v{s}ic \cite{Lv2}. The relation between $V(z)$ in
\eqref{real1} and $W(z)$ in \eqref{een+} is given by
\[
 V(z)=i[W(z)+I]^{-1}[W(z)-I]J.
\]
For an extension of the class of (compactly supported) Herglotz-Nevanlinna
functions in \eqref{nev1} involving a linear term as in
\eqref{nev0}, see \cite{HST1}--\cite{HST3}, \cite{RR1}-\cite{RR2}.
Obviously, general matrix-valued Herglotz-Nevanlinna functions
$V(z)$ cannot be realized in the above mentioned
(Brodski\u{\i}-Liv\v{s}ic) form.

The realization of a different class of Herglotz-Nevanlinna
functions is provided by a linear stationary conservative
dynamical system $\Theta$ of the form
\begin{equation}
\label{col0}
 \Theta =\begin{pmatrix}
          \dA        & K      &  J  \\
          \sH_+\subset\sH\subset\sH_-  &   &\sE
         \end{pmatrix}.
\end{equation}
In this system $\dA$, the \textit{main operator} of the system, is
a bounded linear operator from $\sH_+$ into $\sH_-$ extending a
symmetric (Hermitian) operator $A$ in $\sH$, where
$\sH_+\subset\sH\subset\sH_-$ is a rigged Hilbert space. Moreover,
$K$ is a bounded linear operator from the finite-dimensional
Hilbert space $\sE$ into $\sH_-$, while $J=J^*=J^{-1}$ is acting
on $\sE$, $\varphi_-\in \sE$ is an input vector, $\varphi_+\in
\sE$ is an output vector, and $x\in \sH_+$ is a vector of the
inner state of the system $\Theta$. The system described by
\eqref{col0} is called a {\em canonical} Liv\v{s}ic system or
Brodski\u{\i}-Liv\v{s}ic rigged operator colligation, cf., e.g.
\cite{BT3}-\cite{BHST}. The operator-valued function
\begin{equation}
\label{W1}
 W_\Theta(z)=I-2iK^*(\dA-zI)^{-1}KJ
\end{equation}
is a transfer function (or characteristic function) of the system
$\Theta$. It was shown in \cite{BT3} that a matrix-valued function
$V(z)$ acting on a Hilbert space $\sE$ of the form \eqref{nev0}
can be represented and realized in the form
\begin{equation}
\label{real2}
 V(z)=i[W_\Theta(z)+I]^{-1}[W_\Theta(z)-I]
 =K^*(\dA_{R}-zI)^{-1}K,
\end{equation}
where $W_\Theta(z)$ is a transfer function  of some canonical
scattering ($J=I$) system $\Theta$, and where the ``\textit{real
part}'' $\dA_R=\half(\dA+\dA^*)$ of $\dA$ satisfies $\dA_R \supset
A$ if and only if the function $V(z)$ in \eqref{nev0}
satisfies the following two conditions:
\begin{equation}
\label{cond0} \left\{
\begin{array}{l}
 L = 0, \\
 Qx=\int_{\dR} \frac{t}{1+t^2}\, d\Sigma(t)x
   \quad \mbox{when} \quad
    \int_{\dR} \left(d\Sigma(t)x,x\right)_\sE <\infty.
\end{array}
\right.
\end{equation}
This shows that general matrix-valued Herglotz-Nevanlinna
functions $V(z)$ acting on $\sE$ cannot be realized in the form
\eqref{real2} even by means of a canonical system
(a Brodski\u{\i}-Liv\v{s}ic rigged operator colligation)
$\Theta$ of the form \eqref{col0}.

The main purpose of the present paper is to solve the general
realization problem for matrix-valued Herglotz-Nevanlinna
functions. The case of Herglotz-Nevanlinna functions of the form
\eqref{nev0} with a bounded measure was considered in
\cite{HST1}-\cite{HST3}. In the general case, an appropriate
realization for these functions will be established by introducing
new types systems: so-called {\em non-canonical}
$\Delta_+$-systems and $F_+$-systems. A $\Delta_+$-system or {\em
impedance} system can be written as
\begin{equation}\label{imp1}
    \left\{
\begin{array}{ll}
    (\dD-z F_+)x=K\varphi_-, \\
    \varphi_+=K^*x, \\
\end{array}
\right.
\end{equation}
where $\dD$ and $F_+$ are self-adjoint operators acting from
$\sH_+$ into $\sH_-$ and in addition $F_+$ is an orthogonal
projector in $\sH_+$. In this case the associated transfer
function is given by
\begin{equation}
\label{real03a}
 V(z) =K^*(\dD-z F_+)^{-1}K.
\end{equation}
It will be shown that every matrix-valued Herglotz-Nevanlinna
function  can be represented in the form \eqref{real03a}.

Another type of realization problem deals with so-called
non-canonical \cite{RR1}, \cite{RR2} $F_+$-systems,
\begin{equation}\label{split_system}
    \left\{
\begin{array}{ll}
    (\dA-z F_+)x=KJ\varphi_-, \\
    \varphi_+=\varphi_--2iK^*x, \\
\end{array}
\right.
\end{equation}
also called rigged $F_+$-colligations. This colligation can be
expressed via an array similar to the Brodski\u{\i}-Liv\v{s}ic
rigged operator colligation \eqref{col0}:
\begin{equation}
\label{col1}
 \Theta_{F_+}
         =\begin{pmatrix}
          \dA        &  F_+ & K      &  J  \\
          \sH_+\subset\sH\subset\sH_-  & & &\sE
         \end{pmatrix}.
\end{equation}
The additional ingredient in \eqref{col1} is the operator $F_+$
which is  an
orthogonal projection in $\sH_+$ and $\sH$. The corresponding transfer
function (or $F_+$-characteristic function) is
\begin{equation}
\label{W1+}
 W_{\Theta,F_+}(z)=I-2iK^*(\dA-z F_+)^{-1}KJ.
\end{equation}
It will be shown that every matrix-valued Herglotz-Nevanlinna
function with an invertible matrix $L$ (or $L=0$) in \eqref{nev0}
can be represented in the form
\begin{equation}
\label{real03}
 V(z) =K^*(\dA_R-z F_+)^{-1}K,
\end{equation}
where $\dA_R=\half(\dA+\dA^*)$ is the ``real part'' of the main
operator $\dA$ in the corresponding $F_+$-colligation. The
corresponding $F_+$-characteristic function $W_{\Theta,F_+}(z)$ is
related to the Herglotz-Nevanlinna function $V(z)$ via
\[
 V(z)=i[W_{\Theta,F_+}(z)+I]^{-1}[W_{\Theta,F_+}(z)-I].
\]
Moreover, it will also be shown that the operators $\dD$ and $F_+$
in \eqref{real03a} can be selected so that they satisfy a certain
commutativity condition precisely when the linear term in
\eqref{nev0} is absent, i.e., if $L=0$. When $F_+=I$ the
constructed realization reduces to the Brodski\u{\i}-Liv\v{s}ic
rigged operator colligation (canonical system) \eqref{col0} as
well as to the classical Brodski\u{\i}-Liv\v{s}ic operator
colligation (canonical system) when the measure $d\Sigma(t)$ in
\eqref{nev0} is compactly supported; this includes all the
previous results in the realization problem for matrix-valued
Herglotz-Nevanlinna functions. The results in this paper depend in
an essential way on the theory of extensions in rigged Hilbert
spaces \cite{TSh1}, \cite{Tse}; a concise exposition of this
theory is provided in \cite{TSh1}.

A different approach to realization problems is due to J.A.~Ball
and O.J.~Staffans \cite{BallSt}, \cite{BallSt1},  \cite{St1},
\cite{St2}, \cite{St3}. In particular, they consider canonical
input-state-output systems of the type
\begin{equation}\label{BSsystem}
    \left\{
\begin{array}{ll}
    \dot{x}=Ax(t)+Bu(t), \\
    y(t)=Cx(t)+Du(t), \\
\end{array}
\right.
\end{equation}
with the transfer mapping
\begin{equation}\label{BStransfer}
    T(s)=D+C(sI-A)^{-1}B.
\end{equation}
It follows directly from \cite{BallSt}, \cite{BallSt1}, and
\cite{St3} that for an arbitrary Herglotz-Nevanlinna function
$V(z)$ of the type \eqref{nev0} with $L=0$ the function $-iV(iz)$
can be realized in the form \eqref{BStransfer} by a canonical
impedance conservative system \eqref{BSsystem} considered in
\cite{BallSt}, \cite{BallSt1}, and \cite{St3}. However, this does
not contradict the criteria for the canonical realizations
\eqref{W1}, \eqref{real2}, \eqref{cond0} established by two of the
authors in \cite{BT3} due to the special type ($F_+=I$) of the
Liv\u sic systems (Brodski\u i-Liv\u sic rigged operator
colligations) under consideration. Theorem \ref{thm4} of the
present paper provides a general result for non-canonical
realizations of such functions. The general realization case
involving a non-zero linear term in \eqref{nev0} is also
implicitly treated by Ball and Staffans in \cite{BallSt},
\cite{BallSt1}.

The authors would like to thank Joe Ball and Olof Staffans for
valuable discussions and important remarks.

\section{Some preliminaries}

Let $\sH$ be a Hilbert space with inner product $(x,y)$ and let
$A$ be a closed linear operator in $\sH$ which is Hermitian, i.e.,
$(Ax,y)=(x,Ay)$, for all $x,y\in \dom A$. In general, $A$ need
not be densely defined. The closure of its domain in $\sH$ is
denoted by $\sH_0=\cdom A$. In the sequel $A$ is often considered
as an operator from $\sH_0$ into $\sH$. Then the adjoint $A^*$ of
$A$ is a densely defined operator from $\sH$ into $\sH_0$.
Associated to $A$ are two Hilbert spaces $\sH_+$ and $\sH_-$, the
spaces with a positive and a negative norm. The space $\sH_+$ is
$\dom A^*$ equipped with the graph inner product:
\[
 (f,g)_+=(f,g)+(A^*x,A^*y), \quad f,g \in \dom A^*,
\]
while $\sH_-$ is the corresponding dual space consisting of all
linear functionals on $\sH_+$, which are continuous with respect
to $\|\cdot\|_+$. This gives rise to a triplet $\sH_+\subset
\sH\subset \sH_-$ of Hilbert spaces, which is often called a rigged
Hilbert space associated to $A$. The norms of these spaces satisfy
the inequalities
\[
 \|x\|\leq \|x\|_+, \quad x\in \sH_+,
\quad \mbox{ and } \quad \|x\|_-\leq \|x\|, \quad x\in \sH.
\]
In what follows the prefixes $(+)$-, $(\cdot)$-, and $(-)$- will be
used to refer to corresponding metrics, norms, or inner products of
rigged Hilbert spaces. Recall that there is an isometric operator
$R$, the so-called Riesz-Berezanski\u{\i} operator, which maps
$\sH_-$ onto $\sH_+$ such that
\begin{equation}\label{RBoperator}
\begin{split}
 (x,y)_- & =(x,Ry)=(Rx,y)=(Rx,Ry)_+, \quad x,y \in \sH_-, \\
 (u,v)_+ & =(u,R^{-1}v)=(R^{-1}u,v)=(R^{-1}u,R^{-1}v)_-, \quad u,v \in \sH_+,
\end{split}
\end{equation}
see \cite{Berez}. A closed densely defined linear operator $T$ in
$\sH$ is said to belong to the class $\Omega_A$ if:
\begin{enumerate}
\def\labelenumi{\rm (\roman{enumi})}
\item $A=T\cap T^*$ (i.e. $A$ is the maximal common symmetric part of $T$ and $T^*$);
\item $-i$ is a regular point of $T$.
\end{enumerate}
An operator $\dA\in [\sH_+,\sH_-]$ is called a
$(*)$-\textit{extension} of $T\in \Omega_A$ if the inclusions
\[
 T \subset \dA \quad \mbox{and} \quad T^* \subset \dA^*
\]
are satisfied. Here the adjoints are taken with respect to the
underlying inner products and $[\sH_1,\sH_2]$ stands for the class
of all  linear bounded operators between the Hilbert spaces
$\sH_1$ and $\sH_2$. An operator $\dA\in [\sH_+,\sH_-]$ is a
\textit{bi-extension} of $A$ if $\dA \supset A$ and $\dA^*\supset
A$. Clearly, every $(*)$-extension $\dA$ of $T\in\Omega_A$ is a
bi-extension of $A$.
 A bi-extension $\dA$ of $A$ is called a \textit{self-adjoint bi-extension}
 if $\dA=\dA^*$ and the operator $\wt A$ defined by
\begin{equation}
\label{hatA}
 \wt A =\Big\{\,\{x, \dA x\} :\, x\in \sH_+, \, \dA x\in \sH \,\Big\}.
\end{equation}
is a self-adjoint extension of $A$ in the original Hilbert space $\sH$.
A $(*)$-extension $\dA$ of $T$ is called \textit{correct}
if its ``real part''
$\dA_R:=\half(\dA+\dA^*)$ is a self-adjoint bi-extension of $A$.

For two operators $A$ and $B$ in a Hilbert space $\sH$ the set of
all points $z\in\dC$ such that the operator $(A-zB)^{-1}$ exists
on $\sH$ and is bounded will be denoted by $\rho(A,B)$ and
$\rho(A)=\rho(A,I)$. For some basic facts concerning resolvent
operators of the form $(A-zB)^{-1}$, see \cite{HST3}, \cite{RR1},
\cite{RR2}.

Now proper definitions for both $\Delta_+$-systems and $F_+$-systems
can be given.

\begin{definition}\label{dsystem}
Let $A$ be a closed symmetric operator in a Hilbert space $\sH$
and let $\sH_+\subset \sH \subset \sH_-$ be the rigged Hilbert space
associated with $A$. The system of equations
\begin{equation}\label{imp11}
    \left\{
\begin{array}{ll}
    (\dD-z F_+)x=K\varphi_-, \\
    \varphi_+=K^*x, \\
\end{array}
\right.
\end{equation}
where $\sE$ is a finite-dimensional Hilbert space is called a
$\Delta_+$-system or impedance system if:
\begin{enumerate}
\def\labelenumi{\rm (\roman{enumi})}
\item $\dD\in [\sH_+,\sH_-]$ is a self-adjoint bi-extension of
$A$;
\item $K\in [\sE,\sH_-]$ with $\ker K=\{0\}$ (i.e. $K$ is
invertible);
\item $F_+$ is an orthogonal projection in $\sH_+$ and $\sH$;
\item the set
$\rho(\dD,F_+,K)$ of all points $z\in \dC$ where $(\dD-zF_+)^{-1}$ exists on $\sH\,\cup\,\ran K$
and $(-,\cdot)$-continuous  is open.
\end{enumerate}
\end{definition}

\begin{definition}\label{fsystem}
Let $A$ be a closed symmetric operator in a Hilbert space $\sH$
and let $\sH_+\subset \sH \subset \sH_-$ be the rigged
Hilbert space associated with $A$. The array
\begin{equation}
\label{system1}
 \Theta=\Theta_{F_+}=
         \begin{pmatrix}
          \dA    & F_+  & K & J  \\
          \sH_+\subset \sH \subset \sH_-  & &  & \sE
         \end{pmatrix},
\end{equation}
where $\sE$ is a finite-dimensional Hilbert space
is called an $F_+$-colligation or an $F_+$-system if:
\begin{enumerate}
\def\labelenumi{\rm (\roman{enumi})}
\item $\dA\in [\sH_+,\sH_-]$ is a correct $(*)$-extension of $T\in\Omega_A$;
\item $J=J^*=J^{-1}:\, \sE\to \sE$;
\item $\dA-\dA^*=2i KJK^*$,
where $K\in [\sE,\sH_-]$ and $\ker K=\{0\}$ ($K$ is invertible);
\item $F_+$ is an orthogonal projection in $\sH_+$ and $\sH$;
 \item the set $\rho(\dA,F_+,K)$ of all points $z\in \dC$, where
 $(\dA-zF_+)^{-1}$ exists on $\sH\,\cup\,\ran K$
and $(-,\cdot)$-continuous, is open;
\item the set $\rho(\dA_R,F_+,K)$ of all points $z\in \dC$, where
$(\dA_R-zF_+)^{-1}$ exists on $\sH\,\cup\,\ran K$
and $(-,\cdot)$-continuous,  and the set
      $\rho(\dA,F_+,K) \cap \rho(\dA_R,F_+,K)$ are both open;
\item if $z\in \rho(\dA,F_+,K)$ then $\bar{z}\in
\rho(\dA^*,F_+,K)$; if $z\in \rho(\dA_R,F_+,K)$ then $\bar{z}\in
\rho(\dA_R,F_+,K)$.
\end{enumerate}
\end{definition}
%In the above definitions the property of the type $z\in
%\rho(\dA,F_+,K)$ means that the inverse $(\dA-zF_+)^{-1}$ exists and
%that $\ran K\subset\ran (\dA-zF_+)$, in which case
%$(\dA-zF_+)^{-1}K$ is automatically bounded, since $\sE$ is
%finite-dimensional.

The system \eqref{system1} is conservative in the sense that $\IM
\dA=KJK^*$. It is said to be a scattering system if $J=I$. In this
case the main operator $\dA$ in \eqref{system1} is dissipative:
$\IM \dA\geq 0$. When $F_+=I$ and $\dA$ is a correct
$(*)$-extension of $T\in \Omega_A$ the $F_+$-system in Definition
\ref{fsystem} reduces to a rigged operator colligation (canonical
system) of Brodski\u{\i}-Li\v{v}sic type. It was shown in
\cite{BT3} that each operator $T$ from the class $\Omega_A$ admits
a correct $(*)$-extension $\dA$, which can be included as the main
operator in such a rigged operator colligation and that all the
properties in Definition \ref{fsystem} are automatically
fulfilled. \\

To each $F_+$-system ($F_+$-colligation) in Definition
\ref{fsystem} one can associate a transfer function, or a
characteristic function, via
\begin{equation}
\label{chfunc}
 W_{\Theta}(z)=I-2i K^*(\dA-zF_+)^{-1}KJ.
\end{equation}

\begin{proposition}
\label{Weq} Let $\Theta_{F_+}$ be an $F_+$-colligation of the
form \eqref{system1}. Then for all $z,{w}\in \rho(\dA,F_+,K)$,
\[
\begin{split}
 W_{\Theta_{F_+}}(z)JW_{\Theta_{F_+}}^*(w)-J
 &= 2i(\bar{w}-z)K^*(\dA-zF_+)^{-1}F_+(\dA^*-\bar{w}F_+)^{-1}K, \\
 W_{\Theta_{F_+}}^*(w)JW_{\Theta_{F_+}}(z)-J
 &= 2i(\bar{w}-z)JK^*(\dA^*-\bar{w}F_+)^{-1}F_+(\dA-zF_+)^{-1}KJ.
\end{split}
\]
\end{proposition}

\begin{proof}
By the properties (iii) and (vi) in Definition~\ref{fsystem} one has
for all $z,w\in \rho(\dA,F_+,K)$
\[
\begin{split}
 (\dA-
  & zF_+)^{-1}- (\dA^*-\bar{w}F_+)^{-1} \\
 =& (\dA-zF_+)^{-1}[(\dA^*-\bar{w}F_+)-(\dA-zF_+)]
    (\dA^*-\bar{w}F_+)^{-1}  \\
 =&(z-\bar{w})(\dA-zF_+)^{-1}F_+(\dA^*-\bar{w}F_+)^{-1}
   -2i (\dA-zF_+)^{-1}KJK^*(\dA^*-\bar{w}F_+)^{-1}.
\end{split}
\]
This identity together with \eqref{chfunc} implies that
\[
\begin{split}
 W_{\Theta_{F_+}}(z)
 & J W_{\Theta_{F_+}}^*(w)-J \\
 &=[I-2iK^*(\dA-zF_+)^{-1}KJ]J[I+2iJK^*(\dA^*-\bar{w}F_+)^{-1}K]-J \\
 &=2i(\bar{w}-z)K^*(\dA-zF_+)^{-1}F_+(\dA^*-\bar{w}F_+)^{-1}K.
\end{split}
\]
This proves the first equality. Likewise one proves the second
identity by using
\[
\begin{split}
 (\dA-zF_+)^{-1}-(\dA^*-\bar{w}F_+)^{-1}
  = ( & z-\bar{w})(\dA^*-\bar{w}F_+)^{-1}F_+(\dA^*-zF_+)^{-1} \\
      &-2i (\dA^*-\bar{w}F_+)^{-1}KJK^*(\dA-zF_+)^{-1}.
\end{split}
\]
This completes the proof.
\end{proof}

Proposition~\ref{Weq} shows that the transfer function
$W_{\Theta_{F_+}}(z)$ in \eqref{chfunc} associated to an
$F_+$-system of the form \eqref{system1} is $J$-unitary on
the real axis, $J$-expansive in the upper halfplane, and
$J$-contractive in the lower halfplane with $z\in
\rho(\dA,F_+,K)$. \\

There is another function that one can associate to each
$F_+$-system $\Theta_{F_+}$ of the form \eqref{system1}. It is
defined via
\begin{equation}
\label{Vfunc}
  V_{\Theta_{F_+}}(z)=K^*(\dA_R-zF_+)^{-1}K,
\quad z\in \rho(\dA_R,F_+,K),
\end{equation}
where $\rho(\dA_R,F_+,K)$ is defined above. Clearly,
$\rho(\dA_R,F_+,K)$ is symmetric with respect to the real axis.

\begin{theorem}
\label{Vthm} Let $\Theta_{F_+}$ be an $F_+$-system of the
form \eqref{system1} and let $W_{\Theta_{F_+}}(z)$ and
$V_{\Theta_{F_+}}(z)$ be defined by \eqref{chfunc} and
\eqref{Vfunc}, respectively. Then for all $z,w\in
\rho(\dA_R,F_+,K)$,
\begin{equation}
\label{Viden}
 V_{\Theta_{F_+}}(z)-V_{\Theta_{F_+}}({w})^*
 =(z-\bar{w})K^*(\dA_R-zF_+)^{-1}F_+(\dA_R-\bar{w}F_+)^{-1}K,
\end{equation}
$V_{\Theta_{F_+}}(z)$ is a matrix-valued Herglotz-Nevanlinna function,
and for each $z\in\rho(\dA_R,F_+,K)\cap \rho(\dA,F_+,K)$ the
operators $I+iV_{\Theta_{F_+}}(z)J$ and $I+W_{\Theta_{F_+}}(z)$
are  invertible. Moreover,
\begin{equation}
\label{trans01}
 V_{\Theta_{F_+}}(z)
=i[W_{\Theta_{F_+}}(z)+I]^{-1}[W_{\Theta_{F_+}}(z)-I]J
\end{equation}
and
\begin{equation}
\label{trans02}
 W_{\Theta_{F_+}}(z)
  =[I+iV_{\Theta_{F_+}}(z)J]^{-1}[I-iV_{\Theta_{F_+}}(z)J].
\end{equation}
\end{theorem}

\begin{proof}
For each $z,w\in \rho(\dA_R,F_+,K)$ one has
\begin{equation}
\label{Ares}
 (\dA_R-zF_+)^{-1}-(\dA_R-\bar{w}F_+)^{-1}
  = (z-\bar{w})(\dA_R-zF_+)^{-1}F_+(\dA_R-\bar{w}F_+)^{-1}.
\end{equation}
In view of \eqref{Vfunc} this implies \eqref{Viden}.

%If $\rho(\dA_R,F_+,K)$ is nonempty, then $\cmr\subset
%\rho(\dA_R,F_+,K)$, cf. \cite{HST3}.
Clearly,
\[
 V_{\Theta_{F_+}}(z)^*=V_{\Theta_{F_+}}(\bar{z}).
\]
Moreover, it follows from \eqref{Viden} and the definition
\ref{fsystem} that $V_{\Theta_{F_+}}(z)$  is
%holomorphic on $\rho(\dA_R,F_+,K)$ and has a nonnegative
%imaginary part in the open upper half-plane, so that
%$V_{\Theta_{F_+}}(z)$ is
a matrix-valued Herglotz-Nevanlinna function.

The following identity with
$z\in\rho(\dA,F_+,K)\cap\rho(\dA_R,F_+,K)$
\[
  (\dA_R-zF_+)^{-1}-(\dA-zF_+)^{-1}
  = i(\dA-zF_+)^{-1}\IM \dA (\dA_R-zF_+)^{-1}
\]
leads to
\[
\begin{split}
 K^*(\dA_R-zF_+)^{-1}K-
 & K^*(\dA-zF_+)^{-1}K \\
 & = iK^*(\dA-zF_+)^{-1}KJK^*(\dA_R-zF_+)^{-1}K.
\end{split}
\]
Now in view of \eqref{chfunc} and \eqref{Vfunc}
\[
 2V_{\Theta_{F_+}}(z)+i(I-W_{\Theta_{F_+}}(z))J
   =(I-W_{\Theta_{F_+}}(z))V_{\Theta_{F_+}}(z),
\]
or equivalently, that
\begin{equation}
\label{eq001}
 [I+W_{\Theta_{F_+}}(z)][I+iV_{\Theta_{F_+}}(z)J]=2I.
\end{equation}
Similarly, the identity
\[
  (\dA_R-zF_+)^{-1}-(\dA-zF_+)^{-1}
  = i(\dA_R-zF_+)^{-1}\IM \dA (\dA-zF_+)^{-1}
\]
with $z\in\rho(\dA,F_+,K)\cap\rho(\dA_R,F_+,K)$ leads to
\begin{equation}
\label{eq002}
 [I+iV_{\Theta_{F_+}}(z)J][I+W_{\Theta_{F_+}}(z)]=2I.
\end{equation}
The equalities \eqref{eq001} and \eqref{eq002} show that the
operators are invertible and consequently one obtains
\eqref{trans01} and \eqref{trans02}.
\end{proof}

\section{Impedance Realizations of Herglotz-Nevanlinna functions}
\label{s3}

The realization of Herglotz-Nevanlinna functions has been obtained
for various subclasses. In this section earlier realizations are
combined to present a general realization of an arbitrary
Herglotz-Nevanlinna function by an impedance system. The following
lemma is essentially contained in \cite{HST2};
for completeness a full proof is presented here.
%We present it here to provide further references to the proof.

\begin{lemma}
\label{Dlemma} Let $Q$ be a self-adjoint operator in a
finite-dimensional Hilbert space $\sE$. Then $V(z)=Q$ admits a
representation of the form
\begin{equation}
\label{rep01}
 V(z)=K^*(D-zF_+)^{-1}K, \quad z\in \rho(D,F),
\end{equation}
where $K$ is an invertible mapping from $\sE$ into a Hilbert space
$\sH$, $D$ is a bounded self-adjoint operator in $\sH$, and $F_+$
is an orthogonal projection in $\sH$ whose kernel $\ker F_+$ is
finite-dimensional.
\end{lemma}

\begin{proof}
First assume that $Q$ is invertible.
Let $\sH=\sE$, let $K$ be any invertible mapping
from $\sE$ onto $\sH$, and let $D=KQ^{-1}K^*$. Then $D$ is a
bounded self-adjoint operator in $\sH$. Clearly,
$V(z)=K^*(D-zF)^{-1}K$ with $F=0$, an orthogonal projection in
$\sH$. In the general case, $Q$ can be written as the sum of two
invertible self-adjoint operators $Q=Q^{(1)}+Q^{(2)}$ (for
example, $Q^{(1)}=Q-\varepsilon I$ and $Q^{(2)}=\varepsilon I$,
where $\varepsilon$ is a real number), so that
\[
 Q^{(1)}=K^{(1)*}(D^{(1)}-zF^{(1)})^{-1}K^{(1)},
 \quad
 Q^{(2)}=K^{(2)*}(D^{(2)}-zF^{(2)})^{-1}K^{(2)},
\]
where $K^{(i)}$ is an invertible operator from $\sE$ into a
Hilbert space $\sH^{(i)}=\sE$, $D^{(i)}$ is a bounded self-adjoint
operator in $\sH^{(i)}$, and $F^{(i)}=0$ is an orthogonal
projection in $\sH^{(i)}$, $i=1,2$. (Note that since $K^{(i)}$ is
an arbitrary invertible operator from $\sE$ into $\sH^{(i)}=\sE$
it may as well be chosen as $K^{(i)}=I_{\sE}$). Define
\[
 \sH=\sH^{(1)}\oplus \sH^{(2)},
 \,
 K=\begin{pmatrix} K^{(1)} \\ K^{(2)} \end{pmatrix},
 \,
 D=\begin{pmatrix} D^{(1)} & 0 \\ 0 & D^{(2)} \end{pmatrix},
 \,
 F_+=\begin{pmatrix} F^{(1)} & 0 \\ 0 & F^{(2)} \end{pmatrix}.
\]
Then $K$ is an invertible operator from $\sE$ into the Hilbert
space $\sH$, $D$ is a bounded self-adjoint operator, and $F_+=0$ is
an orthogonal projection in $\sH$. Moreover,
\[
\begin{split}
 Q &= Q^{(1)}+Q^{(2)} \\
   &={K^{(1)}}^*(D^{(1)}-zF^{(1)})^{-1}K^{(1)}
      +{K^{(2)}}^*(D^{(2)}-zF^{(2)})^{-1}K^{(2)} \\
   &= K^*(D-zF_+)^{-1}K,
\end{split}
\]
which proves the lemma.
\end{proof}

Herglotz-Nevanlinna functions of the form \eqref{nev0} which satisfy
the conditions in \eqref{cond0} can be realized by means of the theory of
regularized generalized resolvents, \cite{BT3}, \cite{BT4}. By
means of Lemma \ref{Dlemma} these realizations can be extended to
Herglotz-Nevanlinna functions of the form \eqref{nev0} with $L=0$.

\begin{theorem}
\label{thm3}
Let $V(z)$ be a Herglotz-Nevanlinna function,
acting on a finite-dimensional Hilbert space $\sE$,
with the integral representation
\begin{equation}
\label{nev2}
 V(z)=Q+\int_{\dR}
      \left( \frac{1}{t-z}-\frac{t}{1+t^2}\right)\, d\Sigma(t),
\end{equation}
where $Q=Q^*$ and $\Sigma(t)$ is a nondecreasing matrix-valued
function on $\dR$ satisfying \eqref{int0}.
Then $V(z)$  admits a realization of the form
\begin{equation}
\label{real3}
 V(z)=K^*(\dD-zF_+)^{-1}K,
 \quad z\in \cmr\subset \rho(\dD,F_+,K),
\end{equation}
where $\dD\in [\sH_+,\sH_-]$ is a self-adjoint bi-extension,
$\sH_+\subset\sH\subset \sH_-$ is a rigged Hilbert space, $F_+$ is
an orthogonal projection in $\sH_+$ and $\sH$, $K$ is an injective (invertible)
operator from $\sE$ into $\sH_+$, $K^*\in [\sH_+,\sE]$. Moreover,
the operators $\dD$ and $F_+$ can be selected such that the
following commutativity condition holds:
\begin{equation}
\label{commut1}
 F_-\dD=\dD F_+, \qquad F_-=R^{-1}F_+R \in [\sH_-,\sH_-],
 \end{equation}
 where $R$ is the Riesz-Berezanski\u{\i} operator defined in
 \eqref{RBoperator}.
\end{theorem}

\begin{proof}
According to \cite[Theorem 9]{BT3} each matrix-valued Herglotz-Nevanlinna
function of the form \eqref{nev2} admits a realization of the form
\begin{equation}
\label{rep02}
  V(z)=K^*(\dA_R-zI)^{-1}K=i[W_{\Theta}(z)+I]^{-1}[W_{\Theta}(z)-I],
\end{equation}
where $W_\Theta(z)$ is the transfer function \eqref{W1} of a
system of the form \eqref{col0} if and only if the following
condition holds:
\begin{equation}
\label{cond1a}
  Qx=\int_{\dR} \frac{t}{1+t^2}\, d\Sigma(t)x,
\end{equation}
every vector $x\in\sE$, such that
\begin{equation}
\label{cond1b}
 \int_{\dR} \left(d\Sigma(t)x,x\right)_\sE <\infty.
\end{equation}

To prove the existence of the representation \eqref{real3} for
Herglotz-Nevanlinna functions $V(z)$ which do not satisfy
the condition \eqref{cond1a}, the realization result in Lemma \ref{Dlemma}
will be used.
Denote by $\sE_1$ the linear subspace of vectors
$x\in \sE$ with the property \eqref{cond1b} and let
$\sE_2=\sE\ominus\sE_1$, so that $\sE=\sE_1\oplus\sE_2$. Rewrite
$Q$ in the block matrix form
\[
  Q=\begin{pmatrix} Q_{11} &  Q_{12} \\  Q_{21} & Q_{22} \end{pmatrix},
\quad
  Q_{ij}=P_{\sE_i}Q \uphar {\sE_j}, \quad j=1,2,
\]
and let $\Sigma(t)=(\Sigma_{ij}(t))_{i,j=1}^2$ be decomposed
accordingly. Observe, that by \eqref{int0}, \eqref{cond1a},
\eqref{cond1b} the integrals
\begin{equation}
\label{conv12}
 G_{11}:=\int_{\dR} \frac{t}{1+t^2}\, d\Sigma_{11}(t),
\quad
 G_{12}:=\int_{\dR} \frac{t}{1+t^2}\, d\Sigma_{12}(t)
\end{equation}
are convergent. Let the self-adjoint matrix $G$ be defined by
\begin{equation}
\label{G}
  G=\begin{pmatrix} G_{11} &  G_{12} \\  G_{12}^* & C \end{pmatrix},
\end{equation}
where $C=C^*$ is arbitrary. Now rewrite $V(z)=V_1(z)+V_2(z)$ with
\begin{equation}
\label{V0}
 V_1(z)=Q-G, \quad
 V_2(z)=G+\int_{\dR}
      \left( \frac{1}{t-z}-\frac{t}{1+t^2}\right)\, d\Sigma(t).
\end{equation}
Clearly, for every $x\in \sE_1$ the equality
\[
 Gx=\int_{\dR} \frac{t}{1+t^2}\, d\Sigma(t)x
\]
holds. Consequently, $V_2(z)$ admits the following representation
\begin{equation}
\label{rep03}
  V_2(z)=K_2^*(\dA_R^{(2)}-zI)^{-1}K_2
\end{equation}
where $K_2:\,\sE\to \sH_{-2}$, $K_2^*:\,\sH_{+2}\to \sE$ with
$\sH_{+2}\subset \sH_2 \subset \sH_{-2}$ a rigged Hilbert space,
and where $\dA_R^{(2)}=\half(\dA^{(2)}+(\dA^{(2)})^*)$ is a
self-adjoint bi-extension of a Hermitian operator $A_2$. The
operator $K_2$ is invertible and has the properties
\begin{equation}
\label{Kprop}
\begin{array}{l}
 \ran K_2\subset \ran (\dA^{(2)}-zI),
 \quad
 \ran K_2\subset \ran (\dA_R^{(2)}-zI), \\
 (\dA^{(2)}-zI)^{-1}K_2\in [\sE,\sH_+],
 \quad
 (\dA_R^{(2)}-zI)^{-1}K_2\in [\sE,\sH_+],
\end{array}
\end{equation}
for further details, see \cite{BT3}. Now, by Lemma \ref{Dlemma}
the function $V_1(z)$ admits the representation
\[
 V_1(z)=K_1^*(D_1-zF_{+,1})^{-1}K_1,
\]
where $D_1=D_1^*$ and $F_{+,1}=0$ are acting on a
finite-dimensional Hilbert space $\sH_1=\sE\oplus\sE$ and where
$K_1:\, \sE\to \sH_1$ is invertible. Recall from Lemma
\ref{Dlemma} that
\begin{equation}\label{D1Block}
    D_1=\begin{pmatrix} D_1^{(1)} & 0 \\ 0 & D_1^{(2)}
    \end{pmatrix},\quad
K_1=\begin{pmatrix} K^{(1)}_1 \\ K^{(2)}_1 \end{pmatrix},
\end{equation}
where $K^{(i)}_1:\, \sE\to \sE$, $i=1,2$, and $D_1^{(1)}$,
$D_1^{(2)}$ are defined by means of the decomposition of $Q-G$
into the sum of two invertible self-adjoint operators
$$
Q-G=(Q^{(1)}-G^{(1)})+(Q^{(2)}-G^{(2)}).
$$
Then
\begin{equation}\label{D1i}
    D_1^{(i)}=K_1^{(i)*}(Q^{(i)}-G^{(i)})^{-1}K_1^{(i)},\quad
    i=1,2.
\end{equation}
 To obtain the realization \eqref{real3} for $V(z)$ in
\eqref{nev2}, introduce the following triplet of Hilbert spaces
\begin{equation}\label{RSpace}
\sH_+^{(1)}:=\sE\oplus\sE\oplus \sH_{+2}\subset \sE\oplus\sE\oplus\sH_2
\subset
\sE\oplus\sE\oplus \sH_{-2}:=\sH_-^{(1)},
\end{equation}
i.e., a rigged Hilbert space corresponding to the block
representation of symmetric operator $D_1\oplus A_2$ in
$\sH^{(1)}:=\sH_1\oplus\sH_2$ (where $\sH_1=\sE\oplus\sE$).
Also introduce the following operators
\begin{equation}
\label{KDF}
 \dD=\begin{pmatrix} D_1^{(1)} & 0 &0\\ 0 &D_1^{(2)} & 0\\
  0&0& \dA_R^{(2)} \end{pmatrix},
\quad
 F_+=\begin{pmatrix} 0 & 0 &0\\0& 0 & 0\\0& 0 & I \end{pmatrix},
\quad
 K=\begin{pmatrix} K_1^{(1)}\\K_1^{(2)} \\ K_2 \end{pmatrix}.
\end{equation}
It is straightforward to check that
\begin{equation}
\label{K12}
\begin{split}
 V(z)&=V_1(z)+V_2(z) \\
 &=K_1^{(1)*}(D_1^{(1)}
 -zF_{+,1})^{-1}K_1^{(1)}+K_1^{(2)*}(D_1^{(2)}-zF_{+,1})^{-1}K_1^{(2)}\\
 &\quad +K_2^*(\dA_R^{(2)}-zI)^{-1}K_2 \\
 &= K^*(\dD-zF_+)^{-1}K.
\end{split}
\end{equation}
By the construction,
$A_2 \subset \wt A_R^{(2)}=(\wt A_R^{(2)})^* \subset \dA_R^{(2)}$,
where
\[
 \wt A^{(2)}_R=\{\,\{f,g\}\in \dA_R^{(2)}:\, g\in \sH \,\}
\]
and $A_2$ is a symmetric operator in $\sH_2$, cf. \eqref{hatA}.
Moreover, $\dD$ as
an operator in $[\sH_+^{(1)},\sH_-^{(1)}]$ is self-adjoint,
i.e. $\dD=\dD^*$, and since
\begin{equation}\label{aa2}
 \wh D=\begin{pmatrix} D_1  & 0 \\ 0 & A^{(2)} \end{pmatrix}
 \subset
 \begin{pmatrix} D_1 & 0 \\ 0 & \dA_R^{(2)} \end{pmatrix}
 =\dD,
\end{equation}
and $A =D_1\oplus A_2 \subset \wh D$, the operator $\dD$ is a
self-adjoint bi-extension of the Hermitian operator $A$ in
$\sH_1\oplus\sH_2$. It is easy to see that with  operators in
\eqref{KDF} one obtains the representation \eqref{real3} for
$V(z)$ in \eqref{nev2} and the system constructed with these
operators satisfy the definition \ref{dsystem} of a
$\Delta_+$-system.

Finally, from \eqref{KDF} one obtains $F_-\dD=\dD F_+$, where
$F_+$ and $F_-$ are connected as in \eqref{commut1}. This
completes the proof of the theorem.
\end{proof}

\begin{remark}
According to the recent results by Staffans \cite{St3} an
operator-function $(-i)V(iz)$, where $V(z)$ is defined by
\eqref{nev2} can be realized by an impedance system of the form
\eqref{BSsystem}--\eqref{BStransfer} (see also \cite{BallSt},
\cite{BallSt1}, \cite{St1}, \cite{St2}). This realization is
carried out by using a different approach and does not possess
some of the properties contained in Theorem \ref{thm3}.
\end{remark}

The general impedance realization result for Herglotz-Nevanlinna
functions of the form \eqref{nev0} is now built on Theorem
\ref{thm3} and a representation for linear functions.

\begin{lemma}
\label{linear} Let $L$ be a nonnegative matrix in a
finite-dimensional Hilbert space $\sE$. Then it admits a
realization of the form
\begin{equation}
\label{linear1}
 zL=z {\wh K}^*P {\wh K}=K_3^*(D_3-zF_3)^{-1}K_3,
\end{equation}
where $D_3$ is a self-adjoint matrix in a Hilbert space $\sH_3$,
$P$ is the orthogonal projection onto $\cran L$,
 and $K_3$ is an invertible operator from $\sE$ into $\sH_3$.
\end{lemma}

\begin{proof}
Since $L \geq 0$, there is a unique
nonnegative square root $L^{1/2} \geq 0$ of $L$ with
\[
 \ker L^{1/2}=\ker L, \quad \cran L^{1/2}=\cran L.
\]
Define the operator $\wh K$ in $\sE$ by
\begin{equation}
\label{Khat}
\wh Ku=\left\{
\begin{array}{ll}
    u, & \hbox{$u \in \ker L$;} \\
    L^{1/2}u, & \hbox{$u \in \cran L$.} \\
\end{array}
\right.
 \end{equation}
Then $\wh K$ is invertible and $L^{1/2}=P\wh K$, where $P$
denotes the orthogonal projection onto $\cran L$. Define
\begin{equation}
\label{KDF3}
 \sH_3=\sE\oplus \sE,
 \quad
 K_3=\begin{pmatrix} P \wh K \\ \wh K \end{pmatrix},
 \quad
 D_3=\begin{pmatrix} 0 & iI \\ -iI & 0 \end{pmatrix},
 \quad
 F_{+,3}=\begin{pmatrix} 0 & 0 \\ 0 & I \end{pmatrix}.
\end{equation}
Then $K_3$ is an invertible operator from $\sE$ into $\sH_3$,
$D_3$ is a bounded self-adjoint operator, and $F_{+,3}$ is an
orthogonal projection in $\sH_3$. Moreover,
\begin{equation}
\label{L3rep}
 V_3(z)=zL=z {\wh K}^*P {\wh K}=K_3^*(D_3-zF_{+,3})^{-1}K_3.
\end{equation}
This completes the proof.
\end{proof}

The general realization result for Herglotz-Nevanlinna functions
of the form \eqref{nev0} is now obtained by combining the
earlier realizations.

\begin{theorem}
\label{thm5}
Let $V(z)$ be a matrix-valued Herglotz-Nevanlinna
function in a finite-dimensional Hilbert space $\sE$
with the integral representation
\begin{equation}
\label{V5}
 V(z)=Q+zL+\int_{\dR} \left( \frac{1}{t-z}-\frac{t}{t^2+1} \right) d \Sigma(t),
\end{equation}
where $Q=Q^*$, $L \geq 0$, and $\Sigma(t)$ is a nondecreasing
nonnegative matrix-valued function on $\dR$ satisfying
\eqref{int0}. Then $V(z)$ admits a realization of the form
\begin{equation}
\label{Vrep5}
 V(z)=K^*(\dD-zF_+)^{-1}K
\end{equation}
where $\dD\in [\sH_+,\sH_-]$ is a self-adjoint bi-extension in a
rigged Hilbert space $\sH_+\subset \sH\subset \sH_-$, $F_+$ is an orthogonal
projection in $\sH_+$ and $\sH$, and $K\in [\sE,\sH_-]$ is an invertible
operator from $\sE$ into $\sH_-$.
\end{theorem}

\begin{proof}
Define the following matrix functions
\[
 V_1(z)=Q+\int_{\dR}
            \left( \frac{1}{t-z}-\frac{t}{1+t^2}\right)\, d\Sigma(t),
 \quad
 V_2(z)=zL.
\]
According to Theorem \ref{thm3} the function $V_1(z)$ has a
representation
\[
 V_1(z)=K_1^*(\dD_1-zF_{+,1})^{-1}K_1,
\]
where $\dD_1$, $K_1$ and $F_{+,1}$ are given by the formula
\eqref{KDF}. We recall that $\dD_1$ is a self-adjoint bi-extension
in a rigged Hilbert space $\sH_{-}^{(1)}\subset \sH^{(1)}\subset
\sH^{(1)}_{+}$ given by \eqref{RSpace}, $F_{+,1}$ is an orthogonal
projection in $\sH_+^{(1)}$, and $K_1$ is an invertible mapping
from $\sE$ into $\sH_-^{(1)}$.

According to Lemma \ref{linear} the functions $V_2(z)$ has a
realization of the form \eqref{linear1} with components $\sH_3$,
$D_3$, $K_3$ and $F_{+,3}$ described by \eqref{KDF3}.

Now the final result follows by introducing the rigged Hilbert
space $\sH_3\oplus\sH_{+}^{(1)}\subset\sH_3\oplus\sH^{(1)}
\subset\sH_3\oplus\sH_{-}^{(1)}$ and the operators
\[
 \dD=\begin{pmatrix} D_3 &    0      \\
                   0   & \dD_1
    \end{pmatrix}
   \in [\sH_3\oplus\sH_{+}^{(1)},\sH_3\oplus\sH_{-}^{(1)}],
 \quad
 F_+=\begin{pmatrix} F_{+,3} & 0      \\
                           0   & F_{+,1}
    \end{pmatrix},
 \quad
 K=\begin{pmatrix} K_3 \\
                   K_1
   \end{pmatrix}.
\]
It is straightforward to check that with these operators one
obtains the representation \eqref{Vrep5} for $V(z)$ in \eqref{V5}
and the system constructed with these operators satisfy the
definition \ref{dsystem} of a $\Delta_+$-system.
\end{proof}

For the sake of clarity an extended version for the impedance
realization in the proof of Theorem \ref{thm5} is provided.
The rigged Hilbert space used is
\begin{equation}\label{Greal0}
\sE\oplus\sE\oplus\sE\oplus\sE\oplus \sH_{+2} \subset
\sE\oplus\sE\oplus\sE\oplus\sE\oplus \sH_2 \subset
\sE\oplus\sE\oplus\sE\oplus\sE\oplus \sH_{-2},
\end{equation}
and the operators are given by
\begin{equation}\label{Greal}
   \dD=\begin{pmatrix}
    0&iI&0&0&0\\
    -iI&0& 0&0&0\\
    0& 0 &D_1^{(1)} & 0&0\\
    0& 0 &0 &D_1^{(2)}&0\\
  0&0&0&0& \dA_R^{(2)} \end{pmatrix},\,
 F_+=\begin{pmatrix}
  0&0&0&0&0\\
  0&I&0&0&0\\
  0&0&0&0&0\\
  0&0&0&0&0\\
  0&0&0&0&I
  \end{pmatrix},\,
 K=\begin{pmatrix} P\wh K\\ \wh K \\K_1^{(1)}\\K_1^{(2)}\\ K_2 \end{pmatrix}.
\end{equation}
All the operators in \eqref{Greal} are defined above.

In conclusion of this section it is observed that the general
impedance realization case involving a non-zero linear term in
\eqref{V5} is also implicitly treated by Ball and Staffans in
\cite{BallSt},  \cite{BallSt1}.

\section{$F_+$-system realization results}
\label{s4}

In the general impedance realization results in Theorem \ref{thm3}
and Theorem \ref{thm5} the realizations are in terms of the
operators in \eqref{real3} and \eqref{Vrep5}, respectively. It
remains to identify the Herglotz-Nevanlinna functions as transforms
of transfer functions of appropriate conservative systems.

\begin{theorem}
\label{thm4}
Let $V(z)$ be a Herglotz-Nevanlinna function acting on a
finite-dimensional Hilbert space $\sE$
with the integral representation
\begin{equation}
\label{nev4}
 V(z)=Q+\int_{\dR}
      \left( \frac{1}{t-z}-\frac{t}{1+t^2}\right)\, d\Sigma(t),
\end{equation}
where $Q=Q^*$ and $\Sigma(t)$ is a nondecreasing matrix-valued
function on $\dR$ satisfying \eqref{int0}. Then the function
$V(z)$ can be realized in the form
\begin{equation}
\label{real4}
 V(z)=i[W_{\Theta_{F_+}}(z)+I]^{-1}
  [W_{\Theta_{F_+}}(z)-I],
\end{equation}
where $W_{\Theta_{F_+}}(z)$ is the transfer function given by
\eqref{chfunc} of an $F_+$-system defined in \eqref{system1}. The
$F_+$-system in \eqref{system1} can be taken to be a scattering
system.
\end{theorem}

\begin{proof}
By Theorem \ref{thm3} the function $V(z)$ can be represented in
the form $V(z)=K^*(\dD-zF_+)^{-1}K$, where $K$, $\dD$, and $F_+$
are as in \eqref{KDF} corresponding to the decomposition
\[
 V(z)=V_1(z)+V_2(z),
\]
where
\[
 V_1(z)=Q-G,
 \quad
 V_2(z)=G+\int_{\dR}
      \left( \frac{1}{t-z}-\frac{t}{1+t^2}\right)\, d\Sigma(t),
\]
with a self-adjoint operator $G$ of the form \eqref{G}. With the
notations used in the proof of Theorem \ref{thm3} one may rewrite
$V_1(z)$ and $V_2(z)$ as in \eqref{K12} with
\begin{equation}\label{Depsilon}
    D_1^{(1)}=(Q-G-\varepsilon I)^{-1}, \quad
   D_1^{(2)}=(\varepsilon I)^{-1},\quad
   K_1^{(1)}=\lambda I_\sE, \quad K_1^{(2)}=I_\sE,
\end{equation}
$\dA_R^{(2)}\in [\sH_{+2},\sH_{-2}]$,
$\dA_R^{(2)}=\half(\dA^{(2)}+\dA^{(2)*})$ is associated to a
$(*)$-extension $\dA^{(2)}$ of an operator $T_2\in \Omega_{A_2}$
for which $(-i)\in \rho(T_2)$, cf. \cite{BT3}. The remaining
operators are defined in \eqref{KDF}.

Recall that $K_2$ and the resolvents $(\dA^{(2)}-zI)^{-1}$,
$(\dA_R^{(2)}-zI)^{-1}$ satisfy the properties \eqref{Kprop}. To
construct an $F_+$-system of the form \eqref{system1}
introduce the operator $\dA$ by
\[
  \dA=\dD+iKK^* \in [\sH_+,\sH_-],
\]
where $K$, $\dD$, and $F_+$ are defined in \eqref{KDF}. Then the
block-matrix form of $\dA$ is
\begin{equation}\label{BM_A}
  \dA=\begin{pmatrix} D_1^{(1)}+i\lambda^2I & i\lambda I &i\lambda K_2^*
 \\ i\lambda I &D_1^{(2)}+iI & iK_2^*\\
  i\lambda K_2&iK_2& \dA^{(2)} \end{pmatrix}.
\end{equation}
Let
\begin{equation}
\label{col2}
 \Theta_{F_+} =\begin{pmatrix}
 \dA        &  F_+ & K      &  I  \\
 \sH_+\subset\sH\subset\sH_- & & &\sE
         \end{pmatrix},
\end{equation}
where the rigged Hilbert triplet $\sH_+\subset\sH\subset\sH_-$ is
defined in \eqref{RSpace}, i.e.,
$$\sE\oplus\sE\oplus \sH_{+2}\subset \sE\oplus\sE\oplus\sH_2 \subset
\sE\oplus\sE\oplus \sH_{-2}.$$
 It remains to show that all the properties in Definition
\ref{fsystem} are satisfied. For this purpose, consider the
equation
\[
  (\dA-zF_+)x=(\dD+iKK^*)x-zF_+x=Ke,
\]
or
$$
 \begin{pmatrix} D_1^{(1)}+i\lambda^2I & i\lambda I &i\lambda K_2^*
 \\ i\lambda I &D_1^{(2)}+iI & iK_2^*\\
  i\lambda K_2&iK_2& \dA^{(2)}-zI \end{pmatrix}
  \begin{pmatrix} x_1\\x_2 \\ x_3 \end{pmatrix}=
  \begin{pmatrix} \lambda e\\e \\ K_2e \end{pmatrix}.
$$
Using the decomposition of the operators and taking into account
that
\[
 \dA^{(2)}=\dA_R^{(2)}+iK_2K_2^*
\]
this equation can be rewritten in form of the following system
\begin{equation}
\label{eq10}
\left
 \{
    \begin{array}{l}
     D^{(1)}_1x_1+i\lambda^2Ix_1+i\lambda Ix_2+i\lambda K_2^*x_3=\lambda e, \\
     D^{(2)}_1x_2+i\lambda Ix_1+iIx_2+iK_2^*x_3=e, \\
     (\dA^{(2)}-zI)x_3+i\lambda K_2x_1+iK_2x_2=K_2e.
    \end{array}
\right.
\end{equation}
or
$$
\left
 \{
    \begin{array}{l}
     \frac{1}{\lambda}D^{(1)}_1x_1+i\lambda Ix_1+iIx_2+iK_2^*x_3= e, \\
     D^{(2)}_1x_2+i\lambda Ix_1+iIx_2+iK_2^*x_3=e, \\
     (\dA^{(2)}-zI)x_3+i\lambda K_2x_1+iK_2x_2=K_2e.
    \end{array}
\right.
$$
 In a neighborhood of $(-i)$ the resolvent $(\dA^{(2)}-zI)^{-1}$
is well defined so that by \eqref{Kprop} the third equation in
\eqref{eq10} can be solved for $x_3$:
\begin{equation}
\label{eq11}
  x_3=(\dA^{(2)}-zI)^{-1}K_2e-i(\dA^{(2)}-zI)^{-1}K_2(\lambda x_1+x_2).
\end{equation}
Substitute \eqref{eq11} into the first line of the system yields
$$
\frac{1}{\lambda}D^{(1)}_1x_1+iI(\lambda
x_1+x_2)+K_2^*(\dA^{(2)}-zI)^{-1}K_2(\lambda
x_1+x_2)=e-iK_2^*(\dA^{(2)}-zI)^{-1}K_2e,
$$
Denoting the right hand side by $C$ and using \eqref{chfunc} we
get
$$C=e-iK_2^*(\dA^{(2)}-zI)^{-1}K_2e=\frac{1}{2}\left[I+W_{\Theta_2}(z)\right]e.
$$
Then
$$
\frac{1}{\lambda}D^{(1)}_1x_1+iI(\lambda
x_1+x_2)+K_2^*(\dA^{(2)}-zI)^{-1}K_2(\lambda x_1+x_2)=C,
$$
Multiply both sides by $2i$ and using \eqref{chfunc} one more time
yields
$$
\frac{2i}{\lambda}D^{(1)}_1x_1-\left[I+W_{\Theta_2}(z)\right](\lambda
x_1+x_2)=2iC.
$$
Denoting for further convenience
$B=\left[I+W_{\Theta_2}(z)\right]$ we obtain
$$
\frac{2i}{\lambda}D^{(1)}_1x_1-\lambda Bx_1-Bx_2=2iC.
$$
or
\begin{equation}\label{trying}
  \frac{2i}{\lambda}D^{(1)}_1x_1-\lambda Bx_1-2iC=Bx_2.
\end{equation}
Now we subtract the second equation of the system from the first
and obtain
$$
D^{(1)}_1x_1=\lambda D^{(2)}_1x_2,
$$
or
\begin{equation}\label{x1}
 \lambda(D^{(1)}_1)^{-1}D^{(2)}_1x_2=x_1.
\end{equation}
Applying \eqref{x1} to \eqref{trying}  we get
$$
{2i}D^{(2)}_1x_2-B\lambda^2(D^{(1)}_1)^{-1}D^{(2)}_1x_2-Bx_2=2iC,
$$
and using \eqref{Depsilon}
$$
\frac{2i}{\varepsilon}Ix_2-B[\lambda^2(Q-G-\varepsilon
I)\frac{1}{\varepsilon}+I]x_2=2iC,
$$
or
\begin{equation}\label{Lmatrix}
\Big(2iI-\left[I+W_{\Theta_2}(z)\right][\lambda^2(Q-G)+\varepsilon(1-\lambda^2)
I]\Big)x_2=2i\varepsilon C.
\end{equation}
Choosing $\lambda$ and $\varepsilon$ sufficiently small the matrix
on the left hand side of \eqref{Lmatrix} can be made invertible
for  $z=-i$. Using an invertibility criteria from \cite{Br} we
deduce that \eqref{Lmatrix} is also invertible in a neighborhood
of $(-i)$. Consequently, the system \eqref{eq10} has a unique
solution and $(\dA-zF_+)^{-1}K$ is well defined in a neighborhood
of $(-i)$.

In order to show that the remaining properties in Definition
\ref{fsystem} are satisfied we need to present an operator
$T\in\Omega_A$ such that $\dA$ is a correct $(*)$-extension of
$T$. To construct $T$ we note first that
$(\dA-zF_+)\sH_+\supset\sH$ for some $z$ in a neighborhood of
$(-i)$. This can be confirmed by considering the equation
\begin{equation}\label{Qkernel}
(\dA-zF_+)x=g, \quad x\in\sH_+,
\end{equation}
and showing that it has a unique solution for every $g\in\sH$. The
procedure then is reduced to solving the system \eqref{eq10} with
an arbitrary right hand side $g\in\sH$. Following the steps for
solving \eqref{eq10} we conclude that the system \eqref{Qkernel}
has a unique solution. Similarly one shows that
$(\dA^*-zF_+)\sH_+\supset\sH$. Using the technique developed in
\cite{TSh1} we can conclude that operators $(\dA+iF_+)^{-1}$ and
$(\dA^*-iF_+)^{-1}$ are $(-,\cdot)$-continuous. Define
\begin{equation}\label{DT}
    \begin{array}{l}
    T=\dA, \quad  \dom T=(\dA+iF_+)\sH,\\
    T_1=\dA^*, \quad  \dom T=(\dA^*-iF_+)\sH.
    \end{array}
\end{equation}
One can see that both $\dom T$ and $\dom T_1$ are dense in $\sH$
while operator $T$ is closed in $\sH$. Indeed, assuming that there
is a vector $\phi\in\sH$ that is $(\cdot)$-orthogonal to $\dom T$
and representing $\phi= (\dA^*-iF_+)\psi$ we can immediately get
$\phi=0$. It is also easy to see that $T_1=T^*$. Thus, operator
$T$ defined by \eqref{DT} fits the definition of correct
$(*)$-extension for operator $\dA$. Property (vi) of Definition
\ref{fsystem} follows from Theorem \ref{thm5} and the fact that
$\dA_R=\dD$.

Consequently all the properties for an $F_+$-system $\Theta$
in Definition \ref{fsystem} are fulfilled with the operators and
spaces defined above.
\end{proof}

Now the principal result of the paper will be presented.

\begin{theorem}\label{thm6}
Let $V(z)$ be a matrix-valued Herglotz-Nevanlinna function in a
finite-dimensional Hilbert space $\sE$ with the integral representation
\begin{equation}
\label{V6}
 V(z)=Q+zL+\int_{\dR} \left( \frac{1}{t-z}-\frac{t}{t^2+1}
\right)\, d \Sigma(t),
\end{equation}
where $Q=Q^*$, $L \geq 0$ is an invertible matrix, and $\Sigma(t)$
is a nondecreasing nonnegative matrix-valued function on $\dR$
satisfying \eqref{int0}. Then $V(z)$ can be realized in the form
\begin{equation}
\label{eq20}
 V(z) =i[W_{\Theta_{F_+}}(z)+I]^{-1} [W_{\Theta_{F_+}}(z)-I],
\end{equation}
where $W_{\Theta_{F_+}}(z)$ is a matrix-valued transfer function
of some scattering $F_+$-system of the form \eqref{system1}.
\end{theorem}

\begin{proof}
Decompose the function  $V(z)$ as follows:
\[
 V_1(z)=Q
  + \int_{\dR} \left( \frac{1}{t-z}-\frac{t}{t^2+1} \right) d \Sigma(t)
 \quad\textrm{ and }\quad
 V_2(z)=zL,
\]
and use the earlier realizations for each of these functions.

By Theorem \ref{thm4} the function $V_1(z)$ can be represented by
\[
V_1(z)=i[W_{\Theta_{F_{1,+}}}(z)+I]^{-1}[W_{\Theta_{F_{1,+}}}(z)-I],
\]
where $W_{\Theta_{F_{1,+}}}(z)$ is a matrix-valued transfer
function of some scattering $F_{1,+}$-system,
\begin{equation}
\label{eq26}
 W_{\Theta_{F_{1,+}}}(z)=I-2iK_1^*(\dA_1-zF_{1,+})^{-1}K_1,
\end{equation}
$\dA_1=\dD_1+iK_1K_1^*$ maps $\sH_{+1}$ continuously into
$\sH_{-1}$, $\dD_1$ is a self-adjoint bi-extension, and $\dD_1 \in
[\sH_{+1}, \sH_{-1}]$, $K_1\in [\sE, \sH_{-1}]$.

Following the proof of Theorem \ref{thm5} the function $V_2(z)$
can be represented in the form
\[
 V_2(z)=K_2^*(D_2-zF_{2,+})^{-1}K_2,
\]
where
\begin{equation}
\label{eq28}
  D_2=\begin{pmatrix} 0 & iI \\ -iI & 0 \end{pmatrix}, \quad
  F_{2,+}=\begin{pmatrix} 0 & 0 \\ 0 & I \end{pmatrix},
    \quad
  K_2=\begin{pmatrix} P\widehat{K} \\ \widehat{K} \end{pmatrix},
\end{equation}
and $P$ and $\widehat{K}$ are as in \eqref{Khat}, so that $K_2$ is an
invertible operator from $\sE$ into $\sH_2=\sE \oplus \sE$.
Introduce the triplet
$\sH_{+1}\oplus\sH_2\subset \sH_1\oplus\sH_2\subset
\sH_{-1}\oplus\sH_2$, and consider the operator
\begin{equation}\label{e-28prime}
   \dA=\dD+iK K^*
\end{equation}
 from
$\sH_{+1}\oplus \sH_2$ into $\sH_{-1} \oplus \sH_2$ given by the
block form
\begin{equation}
\label{eq29}
\begin{split}
 \dA
 &=\begin{pmatrix} \dD_1 &0 \\ 0 & D_2 \end{pmatrix}
    +i\begin{pmatrix} K_1 \\ K_2 \end{pmatrix}
   \begin{pmatrix} K_1^* & K_2^* \end{pmatrix} \\
 & =\begin{pmatrix}
       \dA_1 & iK_1K_2^* \\
       iK_2 K_1^* & \dA_2
    \end{pmatrix}.
\end{split}
\end{equation}
Here $\dA_2=D_2+iK_2 K_2^*$.  It will be shown that the equation
\begin{equation}
\label{eq30}
 (\dA-zF_+)x=K e, \quad e \in \sE,
\end{equation}
with
\begin{equation}
\label{eq30-1}
  F_+=\begin{pmatrix} F_{1,+} & 0 \\ 0 & F_{2,+} \end{pmatrix},
  \quad K=\begin{pmatrix} K_1 \\ K_2 \end{pmatrix},
\end{equation}
has always a unique solution $x\in\sH_{+1}\oplus\sH_2$ and
\[
  (\dA-zF_+)^{-1}K \in [\sE, \sH_{+1} \oplus \sH_2].
\]
Taking into account \eqref{eq29}, the equation \eqref{eq30} can be
written as the following system
\begin{equation}
\label{eq31}
  \left\{
  \begin{array}{l}
   (\dA_1-zF_{1,+})x_1+iK_1K_2^*x_2
  =K_1e, \\
  (\dA_2-zF_{2,+})x_2+iK_2K_1^{*}x_1 =K_2e,
  \end{array}
  \right.
\end{equation}
where
\[
 \dA_1=\dD_1+iK_1K_1^*,
 \quad
 \dA_2=D_2+iK_2K_2^*,
\]
By Theorem \ref{thm4} it follows that
\[
  (\dA_1-zF_{1,+})^{-1}K_1 \in [\sE, \sH_{+1}].
\]
Therefore, the first equation in \eqref{eq31} gives
\begin{equation}
\label{eq33}
 x_1=(\dA_1-zF_{1,+})^{-1}K_1e
    -i(\dA_1-zF_{1,+})^{-1}K_1K_2^*x_2.
\end{equation}
Now substituting $x_1$ in the second equation in \eqref{eq31}
yields
\begin{equation}
\label{eq33b}
\begin{split}
 (\dA_2-zF_{2,+})x_2 & +K_2K_1^*(\dA_1-zF_{1,+})^{-1}
      K_1K_2^*x_2  \\
 & =K_2 e
   -iK_2K_1^*(\dA_1-zF_{1,+})^{-1}
    K_1e.
\end{split}
\end{equation}
Taking into account \eqref{eq26}, \eqref{eq28}, and \eqref{eq29}
the identity \eqref{eq33b} leads to
\[
\begin{split}
 &\left(2iI
  -(D_2-zF_+^{(2)})^{-1}K_2
     [I+W_{\Theta_{F_{1,+}}}(z)]K_2^* \right) x_2 \\
 &\qquad =2i(D_2-zF_+^{(2)})^{-1}
   \left( K_2e-iK_2K_1^*
      (\dA_1-zF_+^{(1)})^{-1}K_1e\right).
\end{split}
\]
It will be shown that the matrix-function on the lefthand side,
in front of $x_2$, is invertible. First by straightforward calculations
one obtains
$$
(D_2-zF_+^{(2)})^{-1}=\begin{pmatrix} zI & iI \\ -iI & 0
\end{pmatrix}\in[\sE\oplus\sE,\sE\oplus\sE].
$$
The matrix function $M(z)$ defined by
$$
M(z)=I+W_{\Theta_{F_{1,+}}}(z)\in[\sE,\sE]
$$
is invertible by Theorem~\ref{Vthm}.
It follows from \eqref{eq28} that
$$K_2M(z)=\begin{pmatrix} P\widehat{K} \\ \widehat{K}
\end{pmatrix}M(z)=\begin{pmatrix} L^{1/2} M(z)\\ \widehat{K}M(z)
\end{pmatrix}\in[\sE\oplus\sE,\sE\oplus\sE],
$$
and that
$$
K_2M(z)K_2^*=\begin{pmatrix} L^{1/2} M(z)L^{1/2}&L^{1/2}M(z)\widehat{K}\\
\widehat{K}M(z)L^{1/2}&\widehat{K}M(z)\widehat{K}
\end{pmatrix}\in[\sE\oplus\sE,\sE\oplus\sE].
$$
For any $ 2\times2 $ block-matrix
\[
Z=\begin{pmatrix}
a&b\\c&d\end{pmatrix}
\]
with entries in $[\sE]$ define the matrix-function
$$
N(z)=2iI-(D_2-zF_+^{(2)})^{-1}\begin{pmatrix} a&b\\c&d\end{pmatrix}=
i\begin{pmatrix} 2+zai-c&zbi-d\\a&b+2\end{pmatrix}.
$$
Since the matrix $L>0$ is invertible it follows that $\ker
L=\{0\}$ and $\wh K=L^{1/2}$. Now choose
$$Z=\begin{pmatrix} L^{1/2} M(z)L^{1/2}&L^{1/2}M(z)L^{1/2}\\
L^{1/2}M(z)L^{1/2}&L^{1/2}M(z)L^{1/2}
\end{pmatrix}=\begin{pmatrix}A_0&A_0\\A_0&A_0\end{pmatrix},$$
where $A_0=A_0(z)=L^{1/2} M(z)L^{1/2}$. Note that the matrix-function
$A_0$ is invertible and that $A_0^{-1}=L^{-1/2} M(z)^{-1}L^{-1/2}$.
With this choice of $Z$ one obtains
$$N=N(z)=i\begin{pmatrix} 2I+ziA_0-A_0&ziA_0-A_0\\A_0&A_0+2I\end{pmatrix}.$$
To investigate the invertibility of $N$ consider the system
$$
\begin{pmatrix}
2I+ziA_0-A_0&ziA_0-A_0\\A_0&A_0
+2I\end{pmatrix}\begin{pmatrix}x_1\\x_2\end{pmatrix}
=\begin{pmatrix}0\\0\end{pmatrix},$$ or
\[
\left\{
  \begin{array}{l}
   (2I+ziA_0-A_0)x_1+(ziA_0-A_0)x_2 =0, \\
  A_0x_1+(A_0+2)x_2=0.
  \end{array}
  \right.
\]
Solving the second equation for $x_1$ yields
\[
\left\{
  \begin{array}{l}
   2x_1+ziA_0x_1-A_0x_1+ziA_0x_2-A_0x_2 =0, \\
  x_1=-x_2-2A_0^{-1}x_2.
  \end{array}
  \right.
\]
Substituting $x_1$ into the first equation gives
\[
 (2A_0^{-1}+zi)x_2=0,
\]
%%%%%%%%%%%%%%%%%%%%%%%%%%%%%%%%%%%%%%%%%%%%%%%%%%%%%%%%
or equivalently,
\begin{equation}\label{contra}
 A_0x_2=\frac{2i}{z}x_2.
\end{equation}
Recall that
$$
A_0=A_0(z)=L^{1/2} M(z)L^{1/2}=L^{1/2} [I+
W_{\Theta_1}(z)]L^{1/2}.
$$
For every $z$ in the lower half-plane $W_{\Theta_1}(z)$ is a
contraction (see \cite{BT3}) and thus $\|A_0(z)\|\le 2\|L\|$. This
means that for every $z$ (Im $z<0$) the norm of the left hand side
of \eqref{contra} is bounded while  the norm of the right side can
be made unboundedly large by letting $z\to0$ along the imaginary axis.
This leads to a conclusion that $x_2=0$ and then also $x_1=0$.
Hence, $N=N(z)$ is invertible.
%%%%%%%%%%%%%%%%%%%%%%%%%%%%%%%%%%%%%%%%%%%%%%%%%%%%%%%%

Consequently,
\begin{equation}\label{mmm}
    2iI  -(D_2-zF_+^{(2)})^{-1}K_2[I+W_{\Theta_1}(z)]K_2^*
\end{equation}
is invertible and  $x_2$ depends continuously on $e\in \sE$ in
\eqref{eq33b}, while \eqref{eq33} shows that $x_1$ depends
continuously on $e\in \sE$.

Now we will follow the steps taken in the proof of the Theorem
\ref{thm4} to show that the remaining properties in Definition
\ref{fsystem} are satisfied. We introduce an operator
$T\in\Omega_A$ such that $\dA$ is a correct $(*)$-extension of
$T$.  To construct $T$ we note first that
$(\dA-zF_+)\sH_+\supset\sH$ for some $z$ in a neighborhood of
$(-i)$. This can be confirmed by considering the equation
\begin{equation}\label{Qkernel1}
(\dA-zF_+)x=g, \quad x\in\sH_+,
\end{equation}
and showing that it has a unique solution for every $g\in\sH$. The
procedure then is reduced to solving the system \eqref{eq31} with
an arbitrary right hand side $g\in\sH$. Inspecting the steps of
solving \eqref{eq31} we conclude that the system \eqref{Qkernel1}
has a unique solution. Similarly one shows that
$(\dA^*-zF_+)\sH_+\supset\sH$. Once again relying on \cite{TSh1}
we can conclude that operators $(\dA+iF_+)^{-1}$ and
$(\dA^*-iF_+)^{-1}$ are $(-,\cdot)$-continuous and define
\begin{equation}\label{DT}
    \begin{array}{l}
    T=\dA, \quad  \dom T=(\dA+iF_+)\sH,\\
    T_1=\dA^*, \quad  \dom T=(\dA^*-iF_+)\sH.
    \end{array}
\end{equation}
Using similar to the proof of Theorem \ref{thm4} arguments we note
that both $\dom T$ and $\dom T_1$ are dense in $\sH$ while
operator $T$ is closed in $\sH$.  It is also easy to see that
$T_1=T^*$. Thus, operator $T$ defined by \eqref{DT} fits the
definition of correct $(*)$-extension for operator $\dA$. Property
(vi) of Definition \ref{fsystem} follows from Theorem \ref{thm5}
and the fact that $\dA_R=\dD$.

Therefore, the array
\begin{equation}
\label{col5}
 \Theta_{F_+}
 =\begin{pmatrix}
\dA        & K &F_+     &  I  \\
\sH_{+1}\oplus\sH_2\subset\sH_1\subset\sH_2\subset\sH_{-1}\oplus\sH_2
& &&\sE
  \end{pmatrix}
\end{equation}
is an $F_+$-system and  $V(z)$ admits the realizations
\[
  V(z) =K^*(\dD-zF_+)^{-1}K =i[W_{\Theta_{F_+}}(z)+I]^{-1}
   [W_{\Theta_{F_+}}(z)-I].
\]
This completes the proof.
\end{proof}
It was shown in \cite{HST3} that for the case of compactly
supported measure in \eqref{V6} the function $V(z)$ can be
realized without the restriction on the invertibility of the
linear term $L$.

\section{Minimal Realization}\label{s5}

Recall that a symmetric operator $A$ in a Hilbert space $\sH$ is
called a {\it prime operator} \cite{TSh1}, \cite{BT3} if there
exists no reducing invariant subspace on which it induces a
self-adjoint operator. A notion of a minimal realization is now
defined along the lines of the concept of prime operators.  An
$F_+$-system of the form \eqref{system1} is called
\textit{$F_+$-minimal} if there are no nontrivial reducing invariant
subspaces $\sH^1=\overline{\sH_+^1}$, ($\sH_+^1$ is a $(+)$-subspace
of $\ran F_+$) of $\sH$ where the symmetric operator $A$ induces a
self-adjoint operator. Here the closure is taken with respect to
$(\cdot)$-metric. In the case that $F_+=I$ this definition coincides
with the one used for rigged operator colligations in \cite{Br},
\cite{BT3}.

\begin{theorem}\label{thm5-1}
Let the matrix-valued Herglotz-Nevanlinna function $V(z)$ be
realized in the form
\begin{equation}
\label{eq5-1}
 V(z) =i[W_{\Theta_{F_+}}(z)+I]^{-1}   [W_{\Theta_{F_+}}(z)-I],
\end{equation}
where $W_{\Theta_{F_+}}(z)$ is the transfer function of some
$F_+$-system \eqref{system1}. Then this $F_+$-system can be
reduced to an $F_+$-minimal system of the form \eqref{system1} and
its transfer function gives rise to an $F_+$-minimal realization
of $V(z)$ via \eqref{eq5-1}.
\end{theorem}

\begin{proof}
Let the matrix-valued Herglotz-Nevanlinna function $V(z)$ be
realized in the form \eqref{eq5-1} with an $F_+$-system of the
type \eqref{system1}. Assume that its symmetric operator $A$ has a
reducing invariant subspace $\sH^1=\overline{\sH_+^1}$, ($\sH_+^1$
is a $(+)$-subspace of $\ran F_+$) on which it generates a
self-adjoint operator $A_1$. Then there is the following
$(\cdot,\cdot)$-orthogonal decomposition
\begin{equation}
\label{Adec1} \sH=\sH^0\oplus \sH^1,\quad A=A_0\oplus A_1,
\end{equation}
where $A_0$ is an operator induced by $A$ on $\sH^0$.

The identity \eqref{Adec1} shows that the adjoint of $A$ in $\sH$
admits the orthogonal decomposition $A^*=A_0^*\oplus A_1$. Now
consider  operators $T\supset A$ and  $T^* \supset A$ as in the
definition of the system $\Theta_{F_+}$. It is easy to see that
both $T$ and $T^*$ admit the $(\cdot,\cdot)$-orthogonal
decompositions
\[
T=T_0\oplus A_1,
\]
and
\[
T^\ast=T^\ast_0\oplus A_1,
\]
where $T_0\supset A_0$ and $T^*_0\supset A_0$. Since $T\in
\Omega_A$, the identity $A_0\oplus A_1=T\cap T^*=(T_0\cap
T^\ast_0)\oplus A_1$ holds and $-i$ is a regular point of
$T=T_0\oplus A_1$ or, equivalently, $-i$ is a regular point of
$T_0$. This shows that $T_0\in \Omega_{A_0}$. Clearly,
\[
\sH_+=\sH_+^0 \oplus \sH_+^1=\dom A^\ast_0 \oplus
 \dom A_1.
\]
This decomposition remains valid in the sense of
$(+)$-orthogonality. Indeed, if $f_0\in \sH^0_+$ and
$f_1\in\sH_+^1=\dom A_1 $, then by considering the adjoint of
$A:\sH_0(=\cdom A)\to \sH$ as a mapping from $\sH$ into $\sH_0$
one obtains \[ \aligned
(f_0,f_1)_+&=(f_0,f_1)+(A^\ast f_0,A^\ast f_1)\\
&=(f_0,f_1)+(A_0^\ast f_0,A_1 f_1)\\
&=0+0=0.\\
\endaligned
\]
Consequently, the inclusions $\sH_+\subset\sH\subset\sH_-$ can be
rewritten in the following decomposed forms
\[
\begin{split}%aligned
& \sH_+^0\oplus\sH_+^1 \subset
\sH^0\oplus\sH^1\subset\sH_-^0\oplus\sH_-^1 \\
&=\sH_+^0\oplus\dom A_1
\subset\sH^0\oplus\sH^1\subset\sH_-^0\oplus\sH_-^1.
\end{split}%endaligned
\]
Now let $\dA\in[\sH_+,\sH_-]$ be the correct $(*)$-extension of
$A$ in the definition of the system $\Theta_{F_+}$. Then $\dA$
admits the decomposition $\dA=\dA_0\oplus A_1$ and
$\dA^\ast=\dA_0^*\oplus A_1$. Since $A_1$ is selfadjoint in
$\sH^1$, $\dA_0$ is a correct $(*)$-extension of $T_0$, cf.
\eqref{hatA}. Moreover,
\begin{equation}
\label{Aredu} \aligned \frac{\dA-\dA^\ast}{2i}&=
\frac{(\dA_0\oplus A_1)-(\dA^\ast_0\oplus A_1)} {2i}\\ &=
\frac{\dA_0-\dA_0^\ast}{2i}\oplus\frac{A_1-A_1}{2i}\\ &=
\frac{\dA_0-\dA_0^\ast}{2i}\oplus O,\\
\endaligned
\end{equation}
where $O$ stands for the zero operator. Decompose
$K\in[\sE,\sH_-]$ according to $\sH_-=\sH_-^0\oplus\sH_-^1$ as
follows $K=K_0\oplus K_1$. Then \eqref{Aredu} implies that
\begin{equation}
\label{Kredu} KJK^\ast=K_0JK^\ast_0\oplus O.
\end{equation}
Since $\dim \sE<\infty$ and $\ker K=\{0\}$, one has $\ran K^*=\sE$
and therefore also $\ran JK^*=\sE$. According to \eqref{Kredu}
$K_1(\ran JK^*)=\{0\}$ and therefore $K_1=0$, or equivalently,
$K=K_0\oplus O$. Let $P^0_+$ be the orthogonal projection operator
of $\sH_+$ onto $\sH^0_+$ and let $P^1_+=I-P^0_+$. Then
$K^\ast=K_0^\ast P^0_+$, since for all $f\in \sE$, $g\in \sH_+$ one
has
$$\aligned
(Kf,g)&=(K_0f,g)=(K_0f,g_0+g_1)=(K_0f,g_0)+(K_0f,g_1)\\
&=(K_0f,g_0)=(f,K^\ast_0g_0)=(f,K^\ast_0P^0_+g).\\
\endaligned
$$
Since $\sH_+^1$ is a closed subspace of $\ran F_+$,
$P^0_+=I-P^1_+$ commutes with $F_+$ and therefore
$F_+^0:=F_+P^0_+$ defines an orthogonal projection in $\sH_+^0$.

Now, let $e\in \sE$, let $z\in\rho(\dA,F_+,K)$, and let $x=x^0+x^1
\in \sH_+=\sH_+^0\oplus \sH_+^1$ be such that
$$(\dA-zF_+)x=Ke.$$
Since $K=K_0\oplus O$  the previous identity is equivalent to
$$
(\dA_0\oplus A_1-zF_+)(x^0+x^1)=(K_0\oplus O)e.
$$
Since $F_+x^1=x^1$ and $P^0_+$ commutes with $F_+$, this yields
$$\aligned
(\dA_0-zF_+^0)x^0&=K_0e,\\
(A_1-zI)x^1&=0.
\endaligned$$
It follows from the previous equations that $z\in\rho(A_1)$
because $z\in\rho(\dA,F_+,K)$.  Thus, $\rho(\dA,F_+,K)\subset
\rho(\dA_0,F_+^0,K_0)$ and hence $x^0=(\dA_0-zF_+^0)^{-1}K_0 e$.
On the other hand, $x^0=x= (\dA-zF_+)^{-1}Ke$ and therefore for
all $e\in\sE$ one obtains
$$(\dA-zF_+)^{-1}Ke=(\dA_0-zF_+^0)^{-1}K_0 e$$
and
$$K^\ast (\dA-zF_+)^{-1}Ke=K^\ast_0 (\dA_0-zF_+^0)^{-1}K_0 e.$$
This means that the transfer functions of the system
$\Theta_{F_+}$ in \eqref{system1} and of the system
$$\Theta_{F_+}^0 =\begin{pmatrix} \dA_0 &F_+^0 &K_0 &J\\
\sH_+^0\subset \sH^0 \subset \sH_-^0 &{ }& & \sE
\end{pmatrix}$$ coincide. Therefore, the system $\Theta_{F_+}$
in \eqref{system1} can be reduced to an $F_+$-minimal system of
the same form such that the corresponding transfer functions
coincide. This completes the proof of the theorem.
\end{proof}

The definition of minimality can be extended to $\Delta_+$-systems
in the same manner. Moreover, an $F_+$-system of the form
\eqref{split_system}
$$\left\{
\begin{array}{ll}
    (\dA-z F_+)x=KJ\varphi_-, \\
    \varphi_+=\varphi_--2iK^*x, \\
\end{array}
\right.$$
and a $\Delta_+$-system of the form \eqref{imp1}
$$\left\{
\begin{array}{ll}
    (\dA_R-z F_+)x=K\varphi_-, \\
    \varphi_+=K^*x, \\
\end{array}
\right.$$ where $\dA_R$ is the real part of $\dA$, are minimal (or
non-minimal) simultaneously.

For the $\Delta_+$-systems constructed in Section~\ref{s3} the
minimality can be characterized as follows.

\begin{theorem}\label{cor5-2}
The realization of the matrix-valued Herglotz-Nevanlinna function
$V(z)$ constructed in Theorem~\ref{thm5} is minimal if and only if
the symmetric part $A_2$ of $\dA_R^{(2)}$ defined by \eqref{Greal}
is prime.
\end{theorem}
\begin{proof}
Assume that the system constructed in Theorem \ref{thm5} is not
minimal. Let $\sH^1$ (with $\sH_+^1\subset \ran F_+$) be a
reducing invariant subspace from Theorem \ref{thm5-1} on which $A$
generates a self-adjoint operator $A_1$. Then $\dD=\dD_0\oplus
A_1$ and it follows from the block representations of $\dD$ and
$F_+$ in \eqref{Greal} that $\sH^1$ is necessarily a  subspace of
$\sH_{2}$ in \eqref{Greal0} while $\sH^1_+$ is  a  subspace of
$\sH_{+2}$. To see this let us describe $\ran F_+$ first.
According to \eqref{Greal0}
$\sH_+\subset\sH\subset\sH_-=\sE^4\oplus\sH_{+2}\subset\sE^4\oplus\sH_{2}
\subset\sE^4\oplus\sH_{-2}$ and hence every vector $x\in\sH_+$ can
be written as
$$x=\begin{pmatrix} x_1\\ x_2 \\x_3\\x_4\\ x_5 \end{pmatrix}, \textrm{ where } x_1, x_2,
x_3, x_4\in\sE, x_5\in\sH_{+2}.$$   By \eqref{Greal},
$$F_+x=\begin{pmatrix} 0\\ x_2 \\0\\0\\ x_5 \end{pmatrix},\textrm{ and }
\dD (F_+x)=\begin{pmatrix} ix_2\\ 0 \\0\\0\\ \dA^{(2)}_Rx_5
\end{pmatrix}.$$
This means that $x\in\sH_+^1\subset \ran F_+$ only if $x_2=0$.
Therefore the only possibility for a reducing invariant subspace
$\sH^1$ is to be a subspace of $\sH_2$ while $\sH^1_+$ is a
subspace of $\sH_{+2}$. This proves the claim
$\sH^1_+\subset\sH_{+2}$. Consequently, $\sH^1$ is a reducing
invariant subspace for the symmetric operator $A_2$, in which case
the operator $A_2$ is not prime.

Conversely, if the symmetric operator $A_2$ is not prime, then a
reducing invariant subspace on which $A_2$ generates a
self-adjoint operator is automatically a reducing invariant
subspace for the operator $A$ which belongs to $\ran F_+$. This
completes the proof.
\end{proof}

Finally, Theorem~\ref{cor5-2} implies that a realization of an
arbitrary matrix-valued Herglotz-Nevanlinna function in Theorem
\ref{thm5} can be provided by a minimal $\Delta_+$-system.

\section{Examples}

The paper will be concluded with some simple illustrations of the
main realization result.

\subsection*{Example 1} Consider the following
Herglotz-Nevanlinna function
\begin{equation}\label{e5-1}
    V(z)=1+z-i \tanh \left(\frac{i}{2}z l\right), \quad z \in
\cmr,
\end{equation}
where $l > 0$. An explicit $F_+$-system $\Theta_{F_+}$ will be
constructed so that $V(z)\equiv i[W_{\Theta,F_+}(z)+I]^{-1}[W_{\Theta,F_+}(z)-I]
= V_{\Theta_{F_+}}(z)$. Let the
differential operator $T_2$ in $\sH_2=L^2_{[0,l]}$ be given by
\[
T_2x=\frac{1}{i}\frac{dx}{dt},\quad \dom T_2=\left\{\,x(t) \in
\sH_2:\, x'(t)\in \sH_2, \,x(0)=0 \,\right\},
\]
with adjoint
\[
 T_2^\ast x=\frac{1}{i}\frac{dx}{dt},\quad \dom
T_2^\ast=\left\{\,x(t)\in \sH_2:\, x'(t)\in \sH_2, \,
x(l)=0\,\right\}.
\]
Let $A_2$ be the symmetric operator defined by
\begin{equation}\label{a2}
A_2x=\frac{1}{i}\frac{dx}{dt},\quad \dom A_2=\left\{\,x(t)\in
\sH_2:\, x'(t)\in \sH_2, \,x(0)=x(l)=0 \,\right\},
\end{equation}
with adjoint
$$
A_2^\ast x=\frac{1}{i}\frac{dx}{dt},\quad \dom
A_2^\ast=\left\{\,x(t)\in \sH_2:\, x'(t)\in \sH_2 \,\right\}.
$$
Then $\sH_+=\dom A_2^\ast=W^1_2$ is a Sobolev space with the scalar
product
$$
(x,y)_+=\int^l_0 x(t)\overline{y(t)}\,dt+\int^l_0
x'(t)\overline{y'(t)}\,dt.
$$
Now consider the rigged Hilbert space
$$
W^1_2\subset L^2_{[0,l]}\subset (W_2^1)_-,
$$
and the operators
$$
\aligned
\dA_2x&=\frac{1}{i}\frac{dx}{dt}+i x(0)\left[\delta(x-l)-\delta(x)\right],\\
\dA_2^\ast x&=\frac{1}{i}\frac{dx}{dt}+i x(l)\left[\delta(x-l)-\delta(x)\right],\\
\endaligned
$$
where $x(t)\in W_2^1$ and $\delta(x)$, $\delta(x-l)$ are
delta-functions in $(W^1_2)_-$. Define the operator $K_2$ by
\[
K_2c=c\cdot \frac{1}{\sqrt 2}[\delta(x-l)-\delta(x)], \quad c\in
\dC^1,
\]
so that
\[
K_2^\ast x=\left(x, \frac{1}{\sqrt
2}[\delta(x-l)-\delta(x)]\right)=\frac{1}{\sqrt 2}[x(l)-x(0)],
\]
for $x(t)\in W^1_2$.

Let $D_1=K_1Q_1^{-1}K_1^*=1$, where $Q=1$, and $K_1=1$,
$K_1:\dC\to\dC$. Following \eqref{KDF3} define
$$
\sH_3=\dC\oplus \dC,
 \quad
 K_3=\begin{pmatrix} 1 \\ 1 \end{pmatrix},
 \quad
 D_3=\begin{pmatrix} 0 & iI \\ -iI & 0 \end{pmatrix},
 \quad
 F_{+,3}=\begin{pmatrix} 0 & 0 \\ 0 & I \end{pmatrix}.
$$
Now the corresponding $F_+$-system can be constructed. According
to \eqref{Greal} one has
\[%begin{equation}\label{e5-2}
\dD=\begin{pmatrix}
    0&i&0&0\\
    -i&0& 0&0\\
    0& 0 &1& 0\\
    0&0&0&\dA_{2,R} \end{pmatrix},\,
 F_+=\begin{pmatrix}
  0&0&0&0\\
  0&1&0&0\\
  0&0&0&0\\
  0&0&0&I
  \end{pmatrix},\,
 K=\begin{pmatrix} 1\\ 1\\1\\ K_2 \end{pmatrix},
\]%end{equation}
and it follows from \eqref{e-28prime}-\eqref{eq29} that
\[%begin{equation}\label{e5-3}
 \dA=\dD+iK K^*=\begin{pmatrix}
    i&2i&i&iK_2^*\\
    0&i&i&iK_2^*\\
    i& i &1+i&iK_2^*\\
    iK_2&iK_2&iK_2&\dA_2 \end{pmatrix}.
\]%end{equation}
Consequently, the corresponding $F_+$-system is given by
\begin{equation}\label{e5-4}
\Theta_{F_+}
 =\begin{pmatrix}
\dA        & K      &  F_+&I  \\
\dC^3\oplus W^1_2\subset \dC^3\oplus L^2_{[0,l]}\subset
\dC^3\oplus (W_2^1)_- & &&\dC
  \end{pmatrix},
\end{equation}
where $\dC^3=\dC\oplus \dC\oplus \dC$ and all the operators are
described above. It is well known (see for example \cite{Br}) that
the symmetric operator $A_2$ defined in \eqref{a2} does not have
nontrivial invariant subspaces on which it induces self-adjoint
operators. Thus, the $F_+$-system in \eqref{e5-4} is an
$F_+$-minimal realization of the function $V(z)$ in \eqref{e5-1},
cf. Section \ref{s5}. The transfer function of this system is
$$
W_{\Theta_{F_+}}(z)=\frac{2-i(1+e^{izl})(z+1)}{2e^{izl}+i(1+e^{izl})(z+1)}=
\frac{1-i-zi-\tanh\left(\frac{i}{2}z
l\right)}{1+i+zi+\tanh\left(\frac{i}{2}z l\right)}.
$$

\subsection*{Example 2}
Consider the following Herglotz-Nevanlinna function
\begin{equation}\label{e5-11}
V(z)=\begin{pmatrix} 1 & 0 \\ 0 & 1 \end{pmatrix}+
    z\begin{pmatrix} 1 & 0 \\ 0 & 1 \end{pmatrix}+
    \begin{pmatrix} -i \tanh \left(\pi iz \right) & 0 \\
                     0 & \frac{1-z}{z^2-z-1}
    \end{pmatrix}.
\end{equation}
An explicit $F_+$-system $\Theta_{F_+}$ will be constructed so
that $V(z)\equiv i[W_{\Theta,F_+}(z)+I]^{-1}[W_{\Theta,F_+}(z)-I]=V_{\Theta_{F_+}}(z)$. Let $T_{21}$ be a
differential operator $\sH_2=L^2_{[0,2\pi]}$ given by
$$
T_{21}x=\frac{1}{i}\frac{dx}{dt},\quad \dom
T_{21}=\left\{\,x(t)\in \sH_2:\, x'(t)\in \sH_2,\, x(0)=0
\,\right\},
$$
with adjoint
$$
T_{21}^\ast x=\frac{1}{i}\frac{dx}{dt},\quad \dom
T_{21}^\ast=\left\{\,x(t)\in \sH_2:\, x'(t)\in \sH_2, \, x(2\pi)=0
\,\right\}.
$$
Let $A_{21}$ be the symmetric operator defined by
\begin{equation}\label{a21}
A_{21}x=\frac{1}{i}\frac{dx}{dt},\quad \dom
A_{21}=\left\{\,x(t)\in \sH_2:\, x'(t)\in \sH_2,\,
x(0)=x(2\pi)=0\,\right\},
\end{equation}
with adjoint
$$
A_{21}^\ast x=\frac{1}{i}\frac{dx}{dt},\quad \dom
A_{21}^\ast=\left\{\,x(t)\in \sH_2:\, x'(t)\in \sH_2\,\right\}.
$$
Then $\sH_+=\dom A_{21}^\ast=W^1_2$ is a Sobolev space with the scalar
product
$$(x,y)_+
=\int^{2\pi}_0 x(t)\overline{y(t)}\,dt+\int^{2\pi}_0
x'(t)\overline{y'(t)}\,dt.
$$
Consider the rigged Hilbert space
$$
W^1_2\subset L^2_{[0,2\pi]}\subset (W_2^1)_-,
$$
and the operators
$$
\aligned
\dA_{21}x&=\frac{1}{i}\frac{dx}{dt}+i x(0)\left[\delta(x-2\pi)-\delta(x)\right],\\
\dA_{21}^\ast x&=\frac{1}{i}\frac{dx}{dt}+i x(2\pi)\left[\delta(x-2\pi)-\delta(x)\right],\\
\endaligned
$$
where $x(t)\in W_2^1$ and $\delta(x)$, $\delta(x-2\pi)$ are
delta-functions in $(W^1_2)_-$. Define the operator $K_{21}$ by
\[
K_{21}c=c\cdot \frac{1}{\sqrt 2}[\delta(x-2\pi)-\delta(x)], \quad
c\in \dC^1,
\]
so that
\[
K_{21}^\ast x=\left(x, \frac{1}{\sqrt
2}[\delta(x-2\pi)-\delta(x)]\right)=\frac{1}{\sqrt
2}\,[x(2\pi)-x(0)],
\]
for $x(t)\in W^1_2$. Define
\begin{equation}\label{t22}
T_{22}=\begin{pmatrix}
i&i\\
-i&1
\end{pmatrix}\textrm{ and } K_{22}=\begin{pmatrix} 1 \\
0\end{pmatrix},
\end{equation}
and set
\begin{equation}\label{da2}
\dA_2=\begin{pmatrix}
\dA_{21}&0\\
0&T_{22}
\end{pmatrix}\textrm{ and } K_{2}=\begin{pmatrix} K_{21}&0 \\ 0&K_{22}\end{pmatrix}.
\end{equation}
Now let $D_1=K_1Q^{-1}K_1^*=I_2$, where $Q=I_2$, and $K_1=I_2$,
$K_1:\dC^2\to\dC^2$. Following \eqref{KDF3} define
$$
\sH_3=\dC^4,
 \quad
 K_3=\begin{pmatrix} 1 \\ 1 \\1 \\ 1 \end{pmatrix},
 \quad
 D_3=\begin{pmatrix} 0 &0& i&0 \\ 0 & 0&0&i\\-i &0& 0&0 \\0 &-i& 0&0  \end{pmatrix},
 \quad
 F_{+,3}=\begin{pmatrix} 0&0&0&0\\0&0&0&0\\0&0&1&0\\0&0&0&1  \end{pmatrix}.
$$
Now the corresponding $F_+$-system will be constructed. According
to \eqref{Greal} one has
\begin{equation}\label{e5-2}
\dD=\begin{pmatrix}
    D_3&\vdots&0\\
    \cdots&\cdots& \cdots\\
    0&D_1& 0\\
    \cdots& \cdots &\cdots\\
    0&\vdots&\dA_{2,R} \end{pmatrix},\,
 F_+=\begin{pmatrix}
    F_{+,3}&\vdots&0\\
    \cdots&\cdots& \cdots\\
    0&0& 0\\
    \cdots& \cdots &\cdots\\
    0&\vdots&I \end{pmatrix},\, K=\begin{pmatrix} K_3\\
 \cdots\\K_1\\\cdots\\ K_2 \end{pmatrix},
\end{equation}
and it follows from \eqref{e-28prime}-\eqref{eq29} that
\begin{equation}\label{e5-3}
 \dA=\dD+iK K^*.
\end{equation}
Consequently, the corresponding $F_+$-system is given by
\begin{equation}\label{e5-44}
\Theta_{F_+}
 =\begin{pmatrix}
\dA        & K      &  F_+ &I  \\
\dC^6\oplus W^1_2\subset \dC^6\oplus L^2_{[0,2\pi]}\subset
\dC^6\oplus (W_2^1)_- & &&\dC^2
  \end{pmatrix},
\end{equation}
where all the operators are described above. The transfer function
of this system is given by
\[
W_{\Theta_{F_+}}(z)=
\begin{pmatrix}
\frac{1-i-zi-\tanh(\pi iz)}{1+i+zi+\tanh(\pi iz)}& 0\\
0 & \frac{z^3+iz^2-(3+i)z-i}{-z^3+iz^2+(3-i)z-i}
\end{pmatrix}.
\]
It is easy to see that the maximal symmetric part of the operator
$T_{22}$ in \eqref{t22} is a non-densely defined operator
\begin{equation}\label{a22}
    A_{22}=\begin{pmatrix}
0&i\\
-i&1
\end{pmatrix},
\quad \dom A_{22}=\left\{\,\begin{pmatrix}
0\\
c
\end{pmatrix} :\, c\in \dC \,\right\}.
\end{equation}
 Consequently, the symmetric operator $A_2$ defined by $\dA_2$ in
\eqref{da2}, $\dD$ in \eqref{e5-2}, and $\dA$ in \eqref{e5-3} is
given by
\begin{equation}\label{A2}
\begin{aligned}
   A_2&=\begin{pmatrix}
\frac{1}{i}\frac{dx}{dt}&0&0\\
0&0&i\\
0&-i&1
\end{pmatrix},\\
 \dom A_{2}&=\left\{\,\begin{pmatrix}
x(t)\\
0\\
c
\end{pmatrix} :\, x(t),x'(t)\in \sH_2,\,
x(0)=x(2\pi)=0,\, c\in \dC \,\right\}.
\end{aligned}
\end{equation}
Hence, this operator $A_2$ does not have nontrivial invariant
subspaces on which it induces self-adjoint operators. Thus,
$F_+$-system in \eqref{e5-44} is an $F_+$-minimal realization
of the function $V(z)$ in \eqref{e5-11}, cf. Section \ref{s5}.


\begin{thebibliography}{99}

\bibitem{ADRS}
D.~Alpay, A.~Dijksma, J.~Rovnyak, and H.S.V.~de~Snoo,
\textit{Schur functions, operator colligations, and reproducing
kernel Pontryagin spaces}, Oper. Theory Adv. Appl., 96,
Birkh\"auser Verlag, Basel, 1997.

\bibitem{AlTs1}
D.~Alpay and E.R.~Tsekanovski\u{\i}, ''Interpolation theory in
sectorial Stieltjes classes and explicit system solutions'', Lin.
Alg. Appl., 314 (2000), 91--136.

\bibitem{AlTs2}
D.~Alpay and E.R.~Tsekanovski\u{\i}, ''Subclasses of
Herglotz-Nevanlinna matrix-valued functions and linear systems'',
\textit{Dynamical systems and differential equations}, (ed. J. Du
and S. Hu), An added volume to \textit{Discrete and continuous
dynamical systems}, (2001), 1--14.

\bibitem{A}
Yu.M.~Arlinski\u{\i}, ''On inverse problem of the theory of
characteristic functions of unbounded operator colligations",
Dopovidi Akad. Nauk Ukrain. RSR 2, Ser. A, 105--109.

\bibitem{ArTs0}
Yu.M.~Arlinski\u{\i} and E.R.~Tsekanovski\u{\i}, ''Linear systems
with Schr\"odinger operators and their transfer functions'', Oper.
Theory Adv. Appl., 149, 2004, 47--77.

\bibitem{ArTs1}
Yu.M.~Arlinski\u{\i} and E.R.~Tsekanovski\u{\i}, ''Constant
$J$-unitary factor and operator-valued transfer functions'', {\em
Dynamical Systems and Differential Equations}, (ed. W.Feng, S.Hu
and X.Lu), A supplemental volume to \textit{Discrete and
continuous dynamical systems}, (2003), 48--56.

\bibitem{Aro1}
D.~Arov, ''Passive linear systems and scattering theory'', in
Dynamical Systems, Control Coding, Computer Vision, vol.25 of
Progress in Systems and Control Theory, (1999), Birh\"auser
Verlag, 27--44.

\bibitem{Aro2}
D. Arov, ''Darlington realization of matrix-valued functions'',
Math. USSR Izvestija, 7 (1973), 1295--1326.

\bibitem{AD}
D.~Arov and H.~Dym, "J-inner matrix-functions, interpolation and
inverse problems for canonical systems III.More on the inverse
monodromy problem", Integr. Equat. Oper. Th., 36 (2000), 127--181.

\bibitem{ArGr} D.~Arov and L.Z.~Grossman,
''Scattering matrices in the theory of unitary extensions of
isometric operators'', Math. Nachr., 157, (1992), 105--123.

\bibitem{ArNu} D.~Arov and M.A.~Nudelman,
``Passive linear stationary dynamical scattering systems with
continuous time'', Integral Equat. Oper. Th., 24 (1996), 1--45.

\bibitem{B} J.A.~Ball,
``Linear systems, operator model theory and scattering
multivariable generalizations'', Operator theory: Advances and
Applications, (Winnipeg, MB, 1998), pp. 151--178, Fields Inst.
Commun. 25, Amer. Math. Soc. Providence, RI, 2000.

\bibitem{19} J.A.~Ball and N. Cohen,
``de Branges-Rovnyak operator models and systems theory: a
survey'', Oper. Theory Adv. Appl., 50 (1991), 93--136.

\bibitem{21} J.A.~Ball, I. Gohberg, and L. Rodman,
``Realization and interpolation of rational matrix functions'',
Oper. Theory Adv. Appl., 33 (1988), 1--72.

\bibitem{BGR}
J.A.~Ball, I.~Gohberg, and L.~Rodman, \textit{Interpolation of
rational matrix functions}, Vol. 45, Oper. Theory Adv. Appl.,
Birkh\"auser, 1990.

\bibitem{BallSt}
J.A.~Ball and O.J.~Staffans, ``Conservative state-space
realizations of dissipative system behaviors'', Report No. 37,
Institute Mittag-Leffler, (2002/2003), 55 pp.

\bibitem{BallSt1}
J.A.~Ball and O.J.~Staffans, ``Conservative state-space
realizations of dissipative system behaviors'', Integr. Equ. Oper.
Theory (Online),  Birkh\"auser, 2005, DOI
10.1007/s00020-003-1356-3.

\bibitem{BGK}
H.~Bart, I.~Gohberg, and M.A.~Kaashoek, \textit{Minimal
factorization of matrix and operator functions}, Oper. Theory Adv.
Appl., 1, Birkh\"auser Verlag, Basel, 1979.

\bibitem{BT3}
S.V.~Belyi and E.R.~Tsekanovski\u{\i}, ''Realization theorems for
operator-valued $R$-functions'', Oper. Theory Adv. Appl., 98
(1997), 55--91.

\bibitem{BT4}
S.V.~Belyi and E.R.~Tsekanovski\u{\i}, ''On classes of realizable
operator-valued $R$-functions'', Oper. Theory Adv. Appl., 115
(2000), 85--112.

\bibitem{BHST}
S.V.~Belyi, S.~Hassi, H.S.V. de~Snoo, and E.R.~Tsekanovski\u{\i},
``On the realization of inverse of Stieltjes functions'',
Proceedings of MTNS-2002, University of Notre Dame, CD-ROM, 11p.,
2002.

\bibitem{Berez}
Yu.M.~Berezanski\u{\i}, \textit{Expansion in eigenfunctions of
self-adjoint operators}, vol. 17, Transl. Math. Monographs, AMS,
Providence, 1968.

\bibitem{dBR}
L. de Branges and J.~Rovnyak, ''Canonical models in qunatum
scattering theory'', in \textit{Perturbation theory and its
applications in quantum mechanics}, Wiley \& Sons, New
York-London-Sydney, 1966.

\bibitem{Br}
M.S.~Brodski\u{\i}, \textit{Triangular and Jordan representations
of linear operators}, Moscow, Nauka, 1969 (Russian).

\bibitem{BrLv}
M.S.~Brodski\u{\i} and M.S.~Liv\v{s}ic, ''Spectral analysis of
non-selfadjoint operators and intermediate systems'', Uspekhi
Matem. Nauk, 23, no. 1, 79, (1958), 3--84 (Russian).

\bibitem{Curt4}
R.F.~Curtain and H.~Zwart, \textit{An introduction to
infinite-dimensional linear systems theory}, Springer-Verlag, New
York, 1995.

\bibitem{DoTs}
I.~Dovzhenko and E.R.~Tsekanovski\u{\i}, ''Classes of Stieltjes
operator-functions and their conservative realizations'', Dokl.
Akad. Nauk SSSR, 311 no. 1 (1990), 18--22.

\bibitem{Fur2}
P.A. Fuhrmann, \textit{Linear systems and operators in Hilbert
space}, McGraw-Hill, New York, 1981.

\bibitem{GMKT} F.~Gesztesy, N.J.~Kalton, K.A.~Makarov,
and E.R.~Tsekanovski\u{i}, ''Some applications of operator-valued
Herglotz functions'', Oper. Theory Adv. Appl., 123 (2001),
271--321.

\bibitem{GT}
F.~Gesztesy and E.R.~Tsekanovski\u{\i}, ''On matrix-valued
Herglotz functions'', Math. Nachr., 218 (2000), 61--138.

\bibitem{HST1}
S.~Hassi, H.S.V. de~Snoo, and E.R.~Tsekanovski\u{\i}, ''An
addendum to the multiplication and factorization theorems of
Brodski\u{\i}-Liv\v{s}ic-Potapov'', Applicable Analysis, 77
(2001), 125--133.

\bibitem{HST2}
S.~Hassi, H.S.V. de~Snoo, and E.R.~Tsekanovski\u{\i},
''Commutative and noncommutative representations of matrix-valued
Herglotz-Nevanlinna functions'', Applicable Analysis, 77 (2001),
135--147.

\bibitem{HST3}
S.~Hassi, H.S.V. de~Snoo, and E.R.~Tsekanovski\u{\i},
''Realizations of Herglotz-Nevanlinna functions via
$F$-colligations'', Oper. Theory Adv. Appl., 132 (2002), 183--198.

\bibitem{HST4}
S.~Hassi, H.S.V. de~Snoo, and E.R.~Tsekanovski\u{\i}, ''The
realization problem for Herglotz-Nevanlinna functions'', in
\textit{Unsolved problems in mathematical systems and control
theory}, (ed. V.~Blondel and A.~Megretski), Princeton University
Press, 2004, 8--13.

\bibitem{He}
J.W.~Helton, ''Systems with infinite-dimensional state space: the
Hilbert space approach'', Proc. IEEE, 64 (1976), no. 1, 145--160.

\bibitem{Ku}
M.~Kuijper, \textit{First-order representations of linear
systems}, Birkh\"auser-Verlag, Basel-Boston, 1994.

\bibitem{Lv1}
M.S.~Liv\v{s}ic, ''On spectral decomposition of linear
non-selfadjoint operators'', Math. Sbornik., 30, no.76 (1954),
145--198 (Russian).

\bibitem{Lv2}
M.S.~Liv\v{s}ic, \textit{Operators, oscillations, waves}, Moscow,
Nauka, 1966 (Russian).

\bibitem{LvPo}
M.S.~Liv\v{s}ic and V.P.~Potapov, ''A theorem on the
multiplication of characteristic matrix-functions'', Dokl. Akad.
Nauk SSSR, 72 (1950), 625--628 (Russian).

\bibitem{LiYa}
M.S.~Liv\v{s}ic and A.A.~Yantsevich, \textit{Operator colligations
in Hilbert spaces}, Kharkov University Press, 1971 (Russian)
[English translation: V.H.~Winston \& Sons, 1979]

\bibitem{NF}
B.~Sz.-Nagy and C.~Foias, \textit{Harmonic analysis of operators
on Hilbert space}, North-Holland, New York, 1970.

\bibitem{PoW}
J.W.~Polderman and J.C.~Willems, \textit{Introduction to
mathematical system theory: a behavioural approach}, Springer,
1998.

\bibitem{RR1}
A.~Rutkas and N.~Radbel, ''Linear operator pencils and
noncanonical systems'', Teor. Func. Anal i Prilozhen., 17 (1973),
3--14, (Russian).

\bibitem{RR2}
A.~Rutkas, ``Characteristic function and a model of a linear
operator pencil'', Teor. Func. Func. Anal. i Prilozhen., 45
(1986), 98--111 (Russian) [English Transl., J.Soviet Math., 48
(1990), 451--464].

\bibitem{Sal}
D.~Salamon, ''Infinite dimensional linear systems with unbounded
control and observation: a functional analytic approach'', Trans.
Amer. Math. Soc., 300 (1987), 383--431.

\bibitem{St1}
O.J.~Staffans, \textit{Well-posed linear systems: Part I}, Book
manuscript, available at http://www.abo.fi/$\sim$staffans/, 2002

\bibitem{St2}
O.J.~Staffans, ``Passive and conservative continuous time
impedance and scattering systems, Part I: Well posed systems",
Math. Control Signals Systems, 15, (2002), 291--315.

\bibitem{St3} O.J.~Staffans, Passive and conservative infinite-dimensional
impedance and scattering systems (from a personal point of view),
in Mathematical Systems Theory in Biology, Communication,
Computation, and Finance, J.~Rosenthal and D.~S.~Gilliam, eds, IMA
Volumes in Mathematics and its Applications 134, Springer-Verlag,
New York, 2002, pp. 375-413.

\bibitem{ST3}
O.J.~Staffans and G.~Weiss, ``Transfer functions of regular linear
systems, Part II: the system operator and the Lax-Phillips
semigroup'', Trans. Amer. Math. Soc., 354, (2002), 3229--3262.

\bibitem{Tse}
E.R.~Tsekanovski\u{\i}, ''Accretive extensions and problems on
Stieltjes operator-valued functions realizations'', Oper. Theory
Adv. Appl., 59 (1992), 328--347.

\bibitem{TSh1}
E.R.~Tsekanovski\u{\i} and Yu.L.~Shmul'yan, ''The theory of
bi-extensions of operators on rigged Hilbert spaces. Unbounded
operator colligations and characteristic functions", Russ. Math.
Surv., 32 (1977), 73--131.

\bibitem{WE}
G.~Weiss, ''Transfer functions of regular linear systems. Part I:
characterizations of regularity", Trans. Amer. Math. Soc., 342
(1994), 827--854.

\end{thebibliography}
\end{document}